\documentclass [a4paper,12pt]{article}

\usepackage{amsmath,amsthm}
\usepackage{graphicx, graphics, epsfig,amssymb}
\usepackage[all]{xypic}
\usepackage{floatflt}

\newtheorem{thm}{Theorem}
\newtheorem{lem}{Lemma}
\newtheorem{cor}{Corollary}
\newtheorem{defi}{Definition}
\newtheorem{rem}{Remark}
\newtheorem{exm}{Example}

\newcommand{\R}{\mathbb R}
\newcommand{\N}{\mathbb N}

\DeclareMathOperator{\Hom}{Hom} 
\DeclareMathOperator{\Ob}{Ob}
\DeclareMathOperator{\id}{id} 
\DeclareMathOperator{\const}{const}
 
\DeclareMathOperator{\Sol}{Sol}
\newcommand{\ext}{\mathrm{ext}}
\newcommand{\I}{\mathrm I}
\renewcommand{\d}{\mathrm d}

\newcommand{\PDE}{\mathcal{PDE}}
\newcommand{\PE}{\mathcal{PE}}
\newcommand{\TPE}{\mathcal{TPE}}
\newcommand{\QPE}{\mathcal{QPE}}
\newcommand{\SQPE}{\mathcal{SQPE}}
\newcommand{\AQPE}{\mathcal{AQPE}}
\newcommand{\EPE}{\mathcal{EPE}}
\newcommand{\C}{\mathcal{C}}
\newcommand{\U}{\mathcal{U}}
\newcommand{\Ug}{\U _{\geq}}

\newcommand{\A}{\mathbf{A}}
\newcommand{\B}{\mathbf{B}}

\newcommand{\Abstract}[1]{\vspace{6mm}\par\noindent\hspace*{10mm}
\parbox{140mm}{\small {\bf Abstract.} #1}\par}

\newcommand{\Keywords}[1]{\vspace{3mm}\par\noindent\hspace*{10mm}
\parbox{140mm}{\small {\it Key words:} \rm #1}\par}
\newcommand{\Classification}[1]{\vspace{3mm}\par\noindent\hspace*{10mm}
\parbox{140mm}{\small {\it Mathematics Subject Classification:} \rm #1}\vspace{3mm}\par}

\begin{document}

%================================== Old

\begin{center}
\title{}{\bf
THE STRUCTURE OF THE CATEGORY OF PARABOLIC EQUATIONS
\protect\footnote {The work was partially supported by a RFBR grant 09-01-00139-a (Russia), 
and by the Program for Basic Research of Mathematical Sciences Branch of Russian 
Academy of Sciences. 
This paper was partially written during my stay at IHES; I am grateful to 
this institution for hospitality and for the excellent working conditions. 
I am grateful to Yuri Zarhin and Vladimir Rubtsov for useful remarks.} }
\end{center}

\begin{center}
{\bf Marina Prokhorova
\protect\footnote
{Institute of Mathematics and Mechanics, Russian Academy of Sciences / Ural Branch.
}}

pmf@imm.uran.ru
\end{center}

\Abstract{We define here the category of partial differential equations. 
Special cases of morphisms from an object (equation) are symmetries of the equation 
and reductions of the equation by a symmetry groups, but there are many other morphisms.	
We are mostly interested in a subcategory that arises from 
second order parabolic equations on arbitrary manifolds. 
We introduce a certain structure in this category enabling us to find 
the simplest representative of every quotient object of the given object, 
and develop a special-purpose language for description and study of 
structures of this kind. An example that deals with nonlinear 
reaction-diffusion equation is discussed in more detail.}

\Classification{18B99, 35K05, 35K55, 35K57, 58J70}
%\Classification{18B99, 35K05, 35K10, 35K55, 35K57, 58J65, 58J70.}

\sloppy

\section{Introduction}
\label{sec-intro}

\begin{floatingfigure}{0.25\linewidth}
$\diagram
N \dto_{\pi}|>>\tip \rto^{F} &
N' \dto_{\pi '}|>>\tip
\\
M \rto & M'
\enddiagram$
\caption{}
\label{Fig:Hom0a}

\bigskip

$\diagram
J^k_F(\pi) \dto|<<\tip \drto^{F^k}|>>\tip
\\
J^k(\pi) \dto|>>\tip & 
J^k(\pi') \dto|>>\tip
\\
N \dto_{\pi}|>>\tip \rto^{F}|>>\tip &
N' \dto_{\pi '}|>>\tip
\\
M \rto|>>\tip & M'
\enddiagram$
\caption{}
\label{Fig:Hom0b}
\end{floatingfigure}

We define here the category $\PDE$ of partial differential equations, 
develop a special-purpose language for description and study of its internal structure,
and investigate in detail its full subcategory that arises from 
second order parabolic equations on arbitrary manifolds. 

Let $\pi \colon N \to M$ be a smooth fiber bundle, $E$ be a subset of $k$-jet bundle $J^k(\pi)$.
$E$ could be considered as $k$-th order partial differential equation for a sections of $\pi$, 
that is $s \in \Gamma\pi$ is the solution of $E$ iff its $k$-th prolongation 
$j^k(s) \in \Gamma \left( J^k(\pi) \to M \right)$ is contained in $E$. 
Here $\Gamma\pi$ is the space of smooth sections of $\pi$.

Suppose $\pi \colon N \to M$, $\pi' \colon N' \to M'$ are smooth fiber bundles, 
$F \colon \pi \to \pi'$ is a bundle morphism which is diffeomorphism on the fibers 
(see Fig.~\ref{Fig:Hom0a}).

$F$ induces the map $F^{\ast} \colon \Gamma\pi' \to \Gamma\pi$; 
denote by $\Gamma_F \pi$ its image.
We say that a section of $\pi$ is $F$-projected if it is contained in $\Gamma_F \pi$.
If $F$ is surjective then $F^{\ast}$ is injective,
so we have the map $\Gamma_F\pi \to \Gamma\pi'$.
If $F$ is surjective submersion then we have the map 
$F^k \colon J^k_F(\pi) \to J^k(\pi')$, where $J^k_F(\pi) = F^{\ast}J^k(\pi')$ 
is the bundle of $k$-jets of $F$-projected sections of $\pi$
(see Fig.~\ref{Fig:Hom0b}).

Let us define now the category $\PDE_0$. 
Its objects are pairs $\left( \pi \colon N \to M, E\subset J^k(\pi) \right)$, $k \in \N$,
and morphisms from the object $\left( \pi \colon N \to M, E\subset J^k(\pi) \right)$ 
to the object $\left( \pi' \colon N' \to M', E'\subset J^k(\pi') \right)$ are smooth 
bundle morphisms $F\colon \pi \to \pi'$ satisfying the following conditions:
\begin{enumerate}
	\item $F$ define surjective submersion $M \to M'$,
	\item the diagram Fig.~\ref{Fig:Hom0a} is the pullback square in the category of smooth manifolds, 
	that is for any $x \in M$ the map $\pi^{-1}(x) \to \pi'^{-1}(Fx)$ is diffeomorphism,
	\item $E\cap J^k_F(\pi) = \left(F^k\right)^{-1}(E')$.
\end{enumerate}
If $F \colon (\pi,E) \to (\pi',E')$ is morphism in $\PDE_0$ then $F^{\ast}$ 
define the bijection between the set of all solutions of $E'$ and the set 
of all $F$-projected solutions of $E$.

In Section \ref{sec:catPDE} we define extended category $\PDE$, whose objects are pairs $(N,E)$
where $N$ is a smooth manifold and $E$ is a subset of 
%the total space of 
the bundle $J^k_m(N)$ of $k$-jets of $m$-dimensional submanifolds of $N$, 
and whose morphisms from $(N,E)$ to $(N',E')$ are maps $N \to N'$ 
sutisfying some analogue of conditions (1-3) above.
By definition, the solutions of the equation $E$ are smooth $m$-dimensional 
non-vertical integral manifolds of the Cartan distribution on 
$J^{k}_m(N)$ contained in $E$.
Particularly, the set of solutions includes all $m$-dimensional submanifolds 
$L \subset N$ such that $k$-th prolongation $j^k(L)\subset J^k_m(N)$ 
is contained in $E$.
Any morphism $F \colon (N,E) \to (N',E')$ of $\PDE$ define the bijection 
between the set of all solutions of $E'$ and the set of all 
$F$-projected solutions of $E$ in the same manner as for $\PDE_0$.
%We define here the subcategory $\PDE_0$ of $\PDE$ for the simplicity. 

The category $\PDE$ generalize the notion of a symmetry group in two directions:
\begin{enumerate}
	\item Automorphisms group of the object $(N,E)$ in $\PDE$ is the symmetry group of the equation $E$.
	\item For a symmetry group $G$ of $E$ the natural projection $N \to N/G$ 
	define the morphism $(N,E) \to (N/G,E/G)$ in $\PDE$ where $E/G$ is 
	%the reduced equation that is 
	the equation describing $G$-invariant solutions of $E$.
\end{enumerate}
Note that the morphisms of $\PDE$ go beyond the morphisms of this kind.

Further we want to introduce a certain structure in $\PDE$ 
formed by a lattice of subcategories. 
We receive these subcategories restricting to the equations of specific kind 
(for example, elliptic, parabolic, hyperbolic, linear, quasilinear equations etc.)
or to the morphisms of specific kind (for example, morphisms respecting 
the projection of $N$ on the base manifold $M$ as in $\PDE_0$)
or both.

When we interested in the solutions of some equation 
%from one of these subcategories
it is useful to look for its quotient objects because every quotient object give us the class of solutions of the original equation. 
It may happen that the position of an object in the lattice gives 
information on the morphisms from the object and/or on the kind of the 
simplest representatives of quotient objects. 
In Section \ref{sec:lang} we develop a special-purpose language for description and study of such situations.
We introduce a number of partial orders on the class of all subcategories of fixed category and depict these orders by various arrows (see Table 1 and Fig. 3). 
For instance, we say that a subcategory $\C _1$ is closed in a category $\C $ and depict 
$\xymatrix{{\C} \ar@{^{*}->}[r] & {\C_1}}$
if every morphism in $\C $ with source from $\C _1$ is a morphism in $\C _1$; 
we say that $\C _1$ is plentiful in $\C $ and depict 
$\xymatrix{{\C} \ar@{-->}[r] & {\C_1}}$
if for every $\A  \in \Ob _{\C _1}$ and for
every quotient object of $\A $ in $\C $ there exist a
representative of this quotient object in $\C _1$; and so on. 

We use this language 
%in Section \ref{sec:PE} 
for the detail study of the full subcategory $\PE$ of $\PDE$ that arises from 
second order parabolic equations posed on arbitrary manifolds, 
but we hope that our approach based on category theory may be useful 
for other types of PDE as well. 
Objects of $\PE$ are equations having the form 
$u_t = \sum_{i,j}b^{ij}(t,x,u)u_{ij} + \sum_{i,j}c^{ij}(t,x,u)u_{i}u_{j}
 + \sum_{i} b^{i}(t,x,u)u_{i} + q(t,x,u)$
in a local coordinates $\left(x^i\right)$ on $X$, $X$ is arbitrary smooth manifold.
We prove that every morphism in $\PE$ is of the form
$(t,x,u) \mapsto \left( {{t}'(t),{x'}\left( {t,x}\right),{u}'(t,x,u)} \right)$.
Particularly, $\PE$ appear to be the subcategory of $\PDE_0$.
In Section \ref{sec:subcat} we investigate the structure of $\PE$; 
the part of obtained results is shown on Fig.~\ref{FigBig}, 
where arrows designate specific partial orders on subcategories as was described above. 
On every of three diagrams on Fig.~\ref{FigBig} down movement mean the restriction of class of permitted equations and 
right movement mean the restriction of class of permitted morphisms. 

Using of the developed structure is illustrated in Section 
\ref{sec:reac-diff} on the example of reaction-diffusion equation 
$u_{t} = a(u)\left( {\Delta u + \eta \nabla u} \right) + q(x,u)$, $x \in X$,
posed on a Riemann manifold $X$ equipped with vector field $\eta$; 
let $\A$ be the corresponding object of $\PE$.
If function $a(u)$ is nonlinear enough then  
we receive the following results as immediate corollary of developed investigation of $\PE$ structure:
\begin{itemize}
	\item For every morphism $F\colon \A  \to \A' $ in $\PE$
there exist bijective change of variables in the quotient equation $\A'$
transforming $F$ to the morphism of the form 
$(t,x,u) \to \left( {t,x'(x),u} \right)$ 
and transforming $\A'$ to the equation of the form 
$u'_{t} = a(u')(\Delta u' + \eta'\nabla u') + q'(x',u')$, $x' \in X'$, 
posed on the Riemann manifold $X'$ equipped with vector field $\eta'$, 
with the same function $a$.
	\item If a morphism $F\colon \A  \to \A' $ of the category $\PE$ 
	has the form $(t,x,u) \mapsto \left( {t,{x'}(t,x),{u}'(t,x,u)} \right)$
	and $\A'$ is of the form 
	$u'_{t} = a(u')(\Delta u' + \eta'\nabla u') + q'(x',u')$, 
	then $F$ is of the form $(t,x,u) \to \left( {t,x'(x),u} \right)$.
\end{itemize}

In Sections \ref{sec:catPDE} and \ref{sec:group} we discuss the relations 
between our approach to the factorization of PDE and the other approaches.

Section \ref{sec:Proofs} contains proofs of the theorems.

\section{The category $\PDE_0$ of partial differential equations}
\label{sec:catPDE}

Let $M$, $K$ be smooth manifolds.
A system $E$ of $k$-th order partial differential equations for 
a function $u\colon M \to K$ is given as a system of equations 
$\Phi^l (x,u,\ldots ,u^{(k)})=0$ involving $x$, $u$ and the derivatives 
of $u$ with respect to $x$ up to order $k$, where $x=(x^1,\ldots ,x^m)$ 
are local coordinates on $M$ and $u=(u^1,\ldots ,u^j)$ are local 
coordinates on $K$.
Further we will use a words ``partial differential equation'', ``PDE'' or ``equation'' 
instead of ``a partial differential equation or a system of partial differential equations'' for short.

Remember some things about jets and related notions.
$k$-jet of a smooth function $u \colon M \to K$ at the point $x \in M$
is the equivalence class of smooth functions $M \to K$ 
whose value and partial derivatives up to $k$-th order at $x$ coincide with the ones of $u$.
All $k$-jets of all smooth functions $M \to K$ form the smooth manifold $J^k(M,K)$,
and the natural projection $\pi^k \colon J^{k}(M,K) \to J^{0}(M,K) = M \times K$ 
defines the smooth vector bundle over $M \times K$, which is called $k$-jet bundle.
For every function $u \colon M \to K$ its $k$-th prolongation $j^k(u) \colon M \to J^k(M,K)$ is naturally defined.
$k$-th order PDE for functions $M \to K$ can be considered as subset $E$ of $J^k(M,K)$; 
solutions of $E$ are such functions $u \colon M \to K$ that the image of $j^k(u)$ contained in $E$.

In more general situation we have a smooth fiber bundle $\pi \colon N \to M$ instead of a projection $M \times K \to M$, 
and sections $s \colon M \to N$ instead of functions $u \colon M \to K$. 
By $\Gamma\pi$ denote the space of smooth sections of $\pi$;
remember that a section of $\pi$ is a map $s \colon M \to N$ such that $\pi \circ s$ is the identity.
Definitions of $k$-jet bundle $\pi^k \colon J^{k}(\pi) \to J^{0}(\pi)=N$ and 
of the $k$-th prolongation 
$j^k \colon \Gamma \pi \to \Gamma \left( \pi\circ\pi^k \colon J^k(\pi) \to M \right)$ 
are the same as ones for the functions.
%($\pi$ is a product projection locally over a small regions )
Let $E$ be a subset of $J^k(\pi)$; it could be considered as $k$-th order 
partial differential equation for a sections of $\pi$, 
that is $s \in \Gamma\pi$ is the solution of $E$ if the image of $j^k(s)$ 
is contained in $E$. 

Suppose $\pi \colon N \to M$, $\pi' \colon N' \to M'$ are smooth fiber bundles, 
$F \colon \pi \to \pi'$ is a bundle morphism which is diffeomorphism on the fibers 
(see Fig.~\ref{Fig:Hom0a}).
$F$ induces the map $F^{\ast} \colon \Gamma\pi' \to \Gamma\pi$; 
denote by $\Gamma_F \pi$ its image.
We say that a section of $\pi$ is \textbf{$F$-projected} if it is contained in $\Gamma_F \pi$.
If $F$ is surjective then $F^{\ast}$ is injective,
so we have the map $F_{\#} \colon \Gamma_F\pi \to \Gamma\pi'$.
If $F$ is surjective submersion then we have the map 
$F^k \colon J^k_F(\pi) \to J^k(\pi')$, where $J^k_F(\pi) = F^{\ast}J^k(\pi')$ 
is the bundle of $k$-jets of $F$-projected sections of $\pi$ 
(see Fig.~\ref{Fig:Hom0b}).
Recall that a map $F$ is called a submersion if 
$\d F \colon T_x N \to T_{Fx} N'$ is surjective at each point $x \in N$.

\begin{defi}\label{def_adm_map0}
Let $\pi \colon N \to M$, $\pi' \colon N' \to M'$ be a smooth fiber bundles,
$E$ be a subset of $J^k(\pi)$,
$F \colon \pi \to \pi'$ be a smooth bundle morphism.
We say that $F$ \textbf{is admitted by} $E$ if the following conditions are satisfied:
\begin{enumerate}
	\item $F$ is a surjective submersion,
	\item the diagram Fig.~\ref{Fig:Hom0a} is the pullback square in the category of smooth manifolds, 
	that is for any $x \in M$ the map $\pi^{-1}(x) \to \pi'^{-1}(Fx)$ is diffeomorphism,
	\item $E\cap J^k_F(\pi) = \left(F^k\right)^{-1}(E')$ fore some subset $E'$ of $J^k(\pi')$.
\end{enumerate}
We say that $E'$ is $F$\textbf{-projection} of $E$.
%, as well as a quotient equation for $E$.
\end{defi}

It turns out that the language of category theory is very convenient
for our study of PDE.
Recall that a category $\C $ consists of a collection of
{\it objects} $\Ob _{\C }$, a collection of {\it morphisms}
(or {\it arrows}) $\Hom _{\C }$ and four operations. The
first two operations associate with each morphism $F$ of
$\C $ its {\it source} and its {\it target}, both of which
are objects of $\C $. The remaining two operations are an
operation that associates with each object $\mathbf{C}$ of
$\C $ an {\it identity morphism} ${\rm id}_{\mathbf{C}} \in
\Hom _{\C }$ and an operation of {\it composition} that
associates to any pair $(F,G)$ of morphisms of $\C $ such
that the source of $F$ is coincide with the target of $G$ another
morphism $F \circ G$, their composite. These operations have to
satisfy some natural axioms \cite{Maclain}.

\begin{defi}\label{def:PDE0}
$\PDE_0$ is the category whose objects are pairs 
$\left( \pi \colon N \to M, E\subset J^k(\pi) \right)$, 
$\pi$ is a smooth viber bundle, $k \in \N$,
and morphisms from the object $\left( \pi \colon N \to M, E\subset J^k(\pi) \right)$ 
to the object $\left( \pi' \colon N' \to M', E'\subset J^k(\pi') \right)$ are smooth 
bundle morphisms $F\colon \pi \to \pi'$ admitted by $E$ such that $E'$ is $F$-projection of $E$.
\end{defi}

If $F \colon (\pi,E) \to (\pi',E')$ is morphism in $\PDE_0$ then $F^{\ast}$ 
define the bijection between the set of all solutions of $E'$ and the set 
of all $F$-projected solutions of $E$.

\section{The category $\PDE$ of partial differential equations}
\label{sec:catPDE}

Now we want to define extended category $\PDE$, whose objects are pairs $(N,E)$
where $N$ is a smooth manifold and $E$ is a subset of the bundle $J^k_m(N)$ of 
$k$-jets of $m$-dimensional submanifolds of $N$.

Let $N$ be a $C^r$-smooth manifold, $0<m<\dim N$.
The jet bundle $\pi^k \colon J^{k}_m(N) \to N$ is the fiber bundle 
with the fiber $\left.J^{k}_m(N)\right|_x$ over the point $x \in N$, 
where $\left.J^{k}_m(N)\right|_x$ is the set of equivalence classes 
of smooth $m$-dimensional submanifolds $L$ of $N$ passing through 
$x$ under the equivalence relation of $k$-th order contact in $x$.

$k$-jet of $k$-smooth $m$-dimensional submanifold $L$ over the point 
$x \in L$ is the equivalence class from $\left.J^{k}_m(N)\right|_x$ 
determined by $L$.
Thus we have the prolongation map $j^k\colon L \to J^{k}_m(N)$ 
taking each point $x \in L$ to the $k$-jet of $L$ over $x$ 
(so it is the section of the fiber bundle $J^{k}_m(N)$ restricted to 
$L \subset N$).
For every $k > l \geq 0$ there exist natural projection
$\pi^{k,l}\colon J^{k}_m(N) \to J^{l}_m(N)$
mapping $k$-jet of $L$ to $l$-jet of $L$ over $x$ for every
$m$-dimensional submanifold $L$ of $N$ and every $x \in L$.

For a submanifold $L$ of $N$ the differential of the prolongation 
map $j^k\colon L \to J^{k}_m(N)$ takes the tangent bundle $TL$ 
to the tangent bundle $TJ^{k}_m(N)$.
The closure of the union of the images of $TL$ in $TJ^{k}_m(N)$
when $L$ runs over all $m$-dimensional submanifolds of $N$
is the vector subbundle of $TJ^{k}_m(N)$;
it is called Cartan distribution on $J^{k}_m(N)$.

Let $E$ be a submanifold of $J^k(\pi)$, $\pi \colon N \to M$, $m=\dim N$.
%or a system of such equations for a function $u\colon M \to K$
%(further we will use a word ``system'' for ``a partial differential 
%equation or a system of partial differential equations'' for short).
The graph of a section is $m$-dimensional submanifold of $N$ so 
$J^{k}(\pi)$ is open subspace of $J^{k}_m(N)$ and
$E$ could be considered as the submanifold of $J^{k}_m(N)$. 
The extended version of $E$ is defined as the closure of $E$ 
in $J^{k}_m(N)$ \cite{Olver}.
%A section $u\colon M \to N$ is a solution of $E$ if and only if 
%the $k$-th prolongation of its graph is contained in $\tilde{E}$ .
%This means that a traditional $k$-th order PDE is viewed as a closed submanifold $E$ of $J^k_m(N)$; 
%its solutions are $m$-dimensional submanifolds $L$ of $N$ whose prolongations $j^k(L)$
%lies entirely within $E$
%(we use the same notation both for PDE and for appropriate submanifold 
%of $J^k_m(N)$; the meaning of the notation will be clear from the context).
Because we don't plan to consider infinitesimal properties of $E$ unlike 
the Lie group analysis of PDE, we could consider any subsets $E$ of $J^{k}_m(N)$
as a partial differential equations.
%More generally, arbitrary subset $E$ of $J^{k}_m(N)$ can be 
%considered as a partial differential equation.
By definition, \textbf{solutions} of such an equation are smooth $m$-dimensional 
non-vertical integral manifolds of the Cartan distribution on $J^{k}_m(N)$
contained in $E$.
Note that for any $m$-dimensional submanifold $L$ of $N$
its prolongation $j^k(L)$ is non-vertical integral manifold of the 
Cartan distribution on $J^k_m(N)$.
So if $E \subset J^{k}_m(N)$ is obtained from the traditional PDE 
as it was described above, and $L$ is the graph of a section $u$ of $\pi$,
then $L$ is the solution of $E$ in the above sense if and only if
$u$ is the solution of corresponding traditional PDE in the 
traditional sense.
In addition there is allowed the possibility of both multiply-valued 
solutions and solutions with infinite derivatives
(see \cite{Olver} for the details).
Wherever we write concrete equations in the traditional form below 
we mean the extended versions of these equations that is closure of 
the corresponding set in $J^k_m(N)$.

Now let us introduce some auxiliary notations.

Let $F\colon N \to N'$ be a map. 
We will say that $L \subset N$ is $F$-\textbf{projected} if $L = F^{-1}(F(L))$. 
Note that if $F$ is a surjective submersion and $L$ is $F$-projected 
submanifold of $N$ then $L'=F(L)$ is the submanifold of $N'$.

Let $N$, $N'$ be $C^r$-smooth manifolds, $0<m<\dim N$.
Let $F\colon N \to N'$ be a surjective submersion of smoothness class 
$C^s$, $k \leq s \leq r$.

\begin{defi}\label{def_F-projected}
$F$\textbf{-projected jet bundle} $J_{m,F}^k(N)$ is the submanifold of 
$J^{k}_m(N)$ that consists of $k$-jets of all $m$-dimensional 
$F$-projected submanifolds of $N$.
\end{defi}
We will write $J_{F}^{k}(N)$ instead of $J_{m,F}^{k} \left({N} \right)$ 
if the value of $m$ is clear from the context.

There is natural isomorphism between the bundles $J_{m,F}^k(N)$ and 
$F^{\ast}J_{m'}^k(N')$ over $N$, where 
$F^{\ast}J_{m'}^k(N') = J_{m'}^k(N') \times_{N'} N$ is the pullback of 
$J_{m'}^k(N')$ by $F$, 
$\dim N -m = \dim N' -m'$.
Therefore we can lift the map $F$ to the map
$F^{k}\colon J_{m,F}^{k} \left( {N} \right) \to J^{k}_{m'}\left( {N'} \right)$ 
by the following natural way.
Suppose $\vartheta \in J_{m,F}^{k} \left( {N} \right)$.
\begin{enumerate}
\item
Take an arbitrary $F$-projected submanifold $L$ of $N$ such
that the $k$-prolongation of $L$ pass through $\vartheta$ 
(that is $k$-jet of $L$ over the point $\pi ^k (\vartheta)$ 
is $\vartheta$).
\item
Assign to $\vartheta$ the point 
$\vartheta ' \in J^{k}_{m'} \left( {N'} \right)$, 
where $\vartheta '$ is $k$-jet of the submanifold $L'=F(L)$ of $N'$ 
over the point $F \circ \pi ^{k} \left( {\vartheta}  \right)$.
\end{enumerate}

\begin{figure}[tbh]
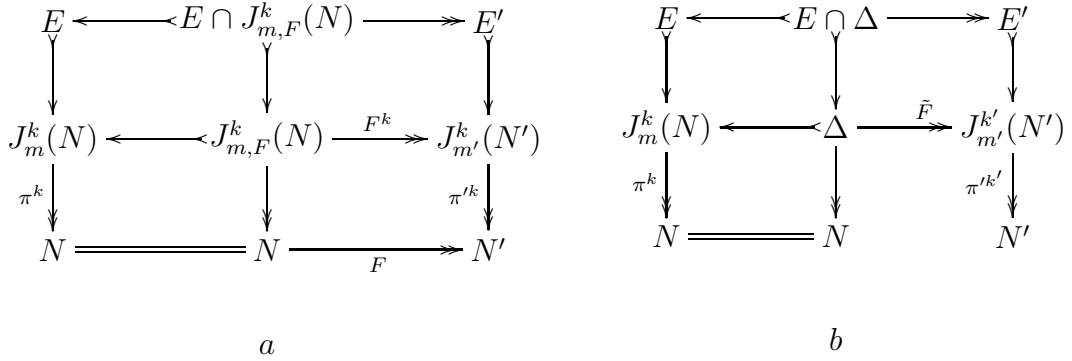

\begin{center}
\[
\diagram
E \dto|<<\tip &
{E \cap J_{m,F}^{k}(N)} \lto|<<\tip \dto|<<\tip \rto|>>\tip &
E' \dto|<<\tip
\\
{J^k_m(N)} \dto_{\pi^k}|>>\tip &
J_{m,F}^{k}(N) \lto|<<\tip \dto|>>\tip \rto^{F^k}|>>\tip &
{J^k_{m'}(N')} \dto_{\pi '^k}|>>\tip
\\
N \ar@2{-}[r]
&
N \rto_{F}|>>\tip & 
N'
\\
& a &
\enddiagram
\quad \quad
\diagram
E \dto|<<\tip &
{E \cap \Delta} \lto|<<\tip \dto|<<\tip \rto|>>\tip &
E' \dto|<<\tip
\\
{J^k_m(N)} \dto_{\pi^k}|>>\tip &
\Delta \lto|<<\tip \dto|>>\tip \rto^{\tilde {F}}|>>\tip &
{J^{k'}_{m'}(N')} \dto_{\pi '^{k'}}|>>\tip
\\
N \ar@2{-}[r]
&
N & 
N'
\\
& b &
\enddiagram
\]
\end{center}
\caption{Morphisms of $\PDE$ (a) and morphisms of $\PDE_{\ext}$ (b)}
\label{FigAdmMap}
\end{figure}

\begin{defi}\label{def_adm_map}
Let $E \subset J^{k}_m(N)$. 
Let $F\colon N \to N'$ be a smooth surjective submersion. 
We say that $F$ \textbf{is admitted by} $E$ if the intersection 
$E \cap J_{m,F}^{k} (N)$ is $F^{k}$-projected subset of $J_{m,F}^{k}(N)$
(see Fig.~\ref{FigAdmMap}(a)). 
Equivalently, $E \cap J_{m,F}^{k}(N)$ is the pre-image 
$\left( {F^{k}} \right)^{-1}\left( {E'} \right)$ 
of some $E'\subset J^{k}_{m'}(N')$; 
we say that $E'$ is $F$\textbf{-projection} of $E$.
\end{defi}

\begin{defi}\label{def_PDE}
A \textbf{category of partial differential equations} $\PDE$ is defined 
as follows:
\begin{itemize}
\item objects of $\PDE$ are pairs $(N,E)$,
where $N$ is a smooth manifold,
$E$ is a subset of $J^{k}_m\left( {N} \right)$ for some integer $k,m \geq 1$;
\item morphisms of $\PDE$ with source $\A =(N,E)$
are all surjective submersions $F\colon N \to N'$ admitted by $E$; 
target of this morphism is $\left( {N',E'} \right)$
where $E'$ is $F$-projection of $E$;
\item the identity morphism from $\A $ is the identity mapping of $N$,
composition of morphisms is composition of appropriate maps.
\end{itemize}
\end{defi}

If $F \colon (N,E) \to (N',E')$ is morphism in $\PDE$ then $F^{\ast}$ 
define the bijection between the set of all solutions of $E'$ and the set 
of all $F^k$-projected solutions of $E$.

In \cite{Prokh1998} we defined the following notion 
of a map admitted by a pair of equations: a map $F\colon N \to N'$
is admitted by an ordered pair of equations $(\A, \A')$, 
$\A  = (N,E)$, $\A' = (N',E')$ if for any $L' \subset N'$
the following two conditions are equivalent:
\begin{itemize}
 \item $L'$ is the graph of a solution of $E'$,
 \item $F^{-1}(L')$ is the graph of a solution of $E$.
\end{itemize}
\noindent
However, we are not happy with this definition; in
particular, because it deals only with global solutions of $E$.
Therefore we now formulate the notion of a map admitted by a equation 
in terms of (locally defined) jet bundles. 

\begin{rem}\label{rem:automor}
Let $\A  = (N,E)$ be an object of $\PDE$. 
Then its automorphism group ${\rm Aut}(\A)$ is the symmetry group 
for the equation $E$.
\end{rem}

\begin{rem}\label{rem:FactByGroup}
Suppose $G$ is a subgroup of the symmetry group of $E$ such that $N/G$ 
is a smooth manifold.
Then the natural projection $N \to N/G$ define the morphism 
$(N,E) \to (N/G,E/G)$ in $\PDE$ where $E/G$ is the equation 
describing $G$-invariant solutions of $E$.
\end{rem}

Therefore, reduction of $E$ by subgroups of ${\rm Aut}(\A)$ 
defines the part of nontrivial morphisms from $\A$.
But the class of all morphisms from $\A$ is significantly richer then the
class of morphisms arising from reduction by subgroups of ${\rm Aut}(\A)$.
Let $\Sol (\A)$ be the set of all solutions of $\A$,
that is of all smooth $m$-dimensional non-vertical integral manifolds
of the Cartan distribution on $J^k_m(N)$ contained in $E$.
In general, the subset
$F^{*}(\Sol(\A '))=\left\{ F^{-1}(L ')\colon L ' \in 
\Sol(\A ') \right\} \subseteq \Sol(\A )$
of solutions of $\A $ arising from a morphism
$F\colon \A \to \A '$ can {\em not\/} be represented as a set of 
solutions that are invariant under some subgroup of ${\rm Aut}(\A)$.
In particular, $F^{*}(\Sol(\A '))$ can be the set of
$G$-invariant solutions, where $G$ is a transformation group that is
not necessarily a symmetry group of $E$. Moreover, for a morphism
$F\colon \A  \to \A '$ it may occur that for every
nontrivial diffeomorphism $g$ of $N$ there exist an element in
$F^{*}(\Sol(\A '))$ that is not $g$-invariant. More detailed
discussion is given in Section \ref{sec:group};
see also \cite{Prokh2000}, \cite{Prokh2001}.

Our approach is conceptually close to the approach developed in
\cite{Elkin} that deals with control systems. If we set aside the
control part and look at this approach relative to ordinary
differential equations, then we get the category of ordinary
differential equations, whose objects are ODE systems of the form
$\dot{x}=\xi$, $x \in X$, where $X$ is a manifold equipped with a
vector field $\xi$, and morphism from a system $\A $ to a
system $\A '$ is a smooth map $F$ from the phase space $X$ of
$\A $ to the phase space $X'$ of $\A '$ that projected
$\xi$ to $\xi'$. In other words, $F$ is a morphism if it transforms
solutions (phase trajectories) of $\A $ to the solutions of
$\A '$: $F_*(\Sol(\A )) = \Sol (\A ')$.

By contrast, we deal with pullbacks of the solutions of the 
quotient equation $\A '$ to the solutions of the original equation $\A $. 
In our approach the number of dependent variables in
the reduced PDE remains the same, while the number of independent0
variables is not increased. Thus in the approach proposed the
quotient object notion is an analogue of the sub-object notion (in
terminology of \cite{Elkin}) with respect to the information about
the solutions of the given equation; however, it is similar to the
quotient object notion with respect to interrelations between the
given and reduced equations.

Note also that described above category of ODE from \cite{Elkin}
is isomorphic to certain subcategory of $\PDE$. Namely, let
us consider the following subcategory $\PDE_1$ of
$\PDE$:
\begin{itemize}
\item objects of $\PDE_1$ are pairs $(N,E)$, where $N=X\times \R$,
$E$ is a first order linear PDE of the form $L_{\xi}u=1$ for unknown
function $u\colon X \to \R$, $\xi \in TX$;
\item morphisms of $\PDE_1$ are morphisms of
$\PDE$ of the form $(x,u)\mapsto (x'(x),u)$.
\end{itemize}
One can easy see that the category of ODE from \cite{Elkin} is
isomorphic to $\PDE_1$: the object $L_{\xi}u=1$ corresponds
to the object $\dot{x}=\xi$, and the morphism $(x,u)\mapsto
(x'(x),u)$ corresponds to the morphism $x\mapsto x'(x)$.

The category of differential equations was also defined in
\cite{Diffeotopy} in a different way: objects are
infinite-dimensional manifolds endowed with integrable
finite-dimensional distribution (particularly, infinite
prolongations of differential equations), and morphisms are smooth
maps such that image of the distribution is contained in the
distribution on the image, similarly to morphisms in \cite{Elkin}.
Thus, the category of differential equations defined in \cite{Diffeotopy} 
is quite different from the category $\PDE$ defined here; 
one should keep it in mind in order to avoid confusion. 
The factorization of PDE $\A $ by a symmetry group described in 
\cite{Diffeotopy} is PDE $\A '$ on the quotient space that described 
images of all solutions of $\A $ at the projection to the quotient space:
$F_*(\Sol(\A )) = \Sol (\A ')$. In that approach every
factorization of $\A$ provides a part of the information about all the
solutions of $\A$. In our approach factorization of $\A$ is such
an equation $\A '$ that the pullbacks of its solutions are
solutions of $\A$: $F^{*}(\Sol(\A ')) \subseteq
\Sol (\A )$; so that from every factorization we obtain the
full information about a certain set of the solutions of the given
equation.

\begin{lem}\label{lem_epi}
All morphisms in $\PDE$ are epimorphisms.
\end{lem}

\begin{figure}[tbh]
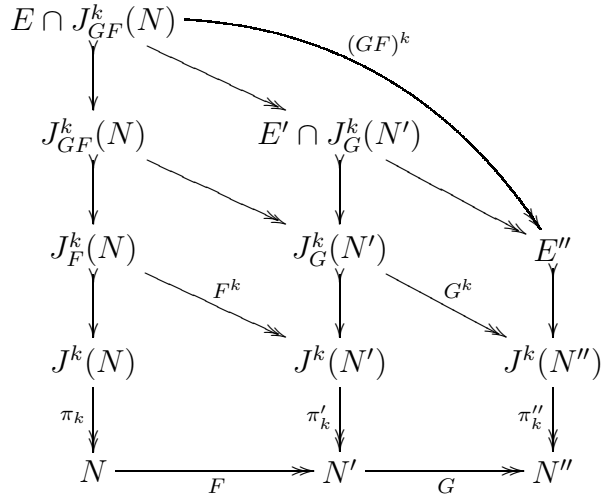

\begin{center}
%\xymatrix{
\[
\diagram
{E \cap J_{GF}^{k}(N)} \dto|<<\tip \drto|>>\tip
\ar@/^2.5pc/[ddrr]^{(GF)^k}|>>\tip
& &
\\
{J_{GF}^{k}(N)} \dto|<<\tip \drto|>>\tip &
{E' \cap J_{G}^{k}(N')} \dto|<<\tip \drto|>>\tip &
\\
{J_{F}^{k}(N)} \dto|<<\tip \drto^{F^k}|>>\tip &
{J_{G}^{k}(N')} \dto|<<\tip \drto^{G^k}|>>\tip &
{E''} \dto|<<\tip
\\
{J^k(N)} \dto_{\pi_k}|>>\tip &
{J^k(N')} \dto_{\pi '_k}|>>\tip &
{J^k(N'')} \dto_{\pi ''_k}|>>\tip
\\
{N} \rto_{F}|>>\tip &
{N'} \rto_{G}|>>\tip &
{N''}
\enddiagram
\]
%}
\caption{Diagram for Lemma \ref{lem_composition}}
\label{FigCompAdm}
\end{center}
\end{figure}

\begin{lem}\label{lem_composition}
Suppose $(N,E)$, $(N',E')$, $(N'',E'')$ are objects of $\PDE$,
$F\colon N \to N'$ is a morphism from $(N,E)$ to $(N',E')$ in $\PDE$, 
$G\colon N' \to N''$ is surjective submersion. 
Then the following two conditions are equivalent 
(see Fig.~\ref{FigCompAdm}):
\begin{itemize}
  \item $G$ is a morphism from $(N',E')$ to $(N'',E'')$,
  \item $GF$ is a morphism from $(N,E)$ to $(N'',E'')$.
\end{itemize}
\end{lem}

\section{The extended category $\PDE_{\ext}$ of partial differential equations}

Note that the Cartan distribution $C^k(N)$ on $J^{k}_m\left( {N} \right)$
restricted to $J_{m,F}^{k} \left( {N} \right)$ coincides with the
lifting $\left( {F^{k}} \right)^{\ast} C_{m'}^{k}\left( {N'} \right)$ of
the Cartan distribution on $J_{m'}^{k}\left( {N'} \right)$, $m'=m-\dim N+\dim N'$. 
Taking this into account and using the analogy with higher symmetry group, 
we replace $J^k_{m,F}(N)$ to arbitrary submanifold $\Delta$ of $J_m^{k}(N)$. 
Thus we obtain the category $\PDE_{\ext}$ with the same objects as $\PDE$ and 
extended set of morphisms involving transformations of jets. 
(This category will not be used in the rest of the paper.)

\begin{defi}\label{def:PDEext}
An \textbf{extended category of partial differential equations} $\PDE_{\ext}$ 
is defined as follows:
\begin{itemize}
\item objects of $\PDE_{\ext}$ are pairs $(N,E)$, where $N$ is a smooth manifold,
$E$ is a subset of $J^{k}_m\left( {N} \right)$ for some integer $k,m \geq 1$;
\item morphisms of $\PDE_{\ext}$ from $\A =(N,E \subset J_{m}^{k}(N))$ to 
$\A' =(N',E' \subset J_{m'}^{k'}(N'))$
are all pairs $\left( \Delta, \tilde {F} \right)$ such that 
$\Delta$ is a smooth submanifold of $J_m^{k}(N)$,
$\tilde {F}\colon\Delta \to J_{m'}^{k'}(N')$ is a surjective submersion,
the Cartan distribution on $J_m^{k}\left( {N} \right)$ restricted to $\Delta $
coincide with the lifting $\tilde {F}^{\ast} C_{m'}^{k'}\left( {N'} \right)$
of the Cartan distribution on $J_{m'}^{k'}\left( {N'} \right)$,
and $E \cap \Delta = \tilde{F}^{-1}(E')$
(see Fig.~\ref{FigAdmMap}(b))
\item the identity morphism from $\A $ is $\left( {\Delta=J^k_m(N), \tilde {F}=\id_N} \right)$, 
composition of $\left( {\Delta\subset J_m^{k}(N), \tilde{F}\colon\Delta \to J_{m'}^{k'}(N')} \right)$ 
and $\left( {\Delta'\subset J_m'^{k'}(N'), \tilde{F}'\colon\Delta'\to J_{m''}^{k''}(N'')} \right)$ 
is $\left( {\tilde{F}^{-1}(\Delta'), \tilde{F'}\circ\tilde{F}} \right)$.
\end{itemize}
\end{defi}

For each integral manifold of the Cartan distribution on $E'$ its
preimage is an integral manifold of the Cartan distribution on $E$,
so for each solution of $E'$ its pullback is some solution of $E$.

$\PDE$ embeds to $\PDE_{\ext}$ by the following natural way:
to the morphisms $F \colon N \to N'$ of $\PDE$ from the equation of $k$-th order 
we assign the morphisms $\left( {\Delta, \tilde {F}} \right)$ of $\PDE_{\ext}$ 
such that $\Delta=J_{m,F}^{k}(N)$, $\tilde {F} = F^k$.

\section{Usage of subcategories}
\label{sec:lang}

We start with review of some basic definitions of category theory
\cite{Maclain}. Given a category $\C $ and an object
$\A $ of $\C $, one may construct the category
$\left( {\A  \downarrow \C } \right)$ of objects
under $\A $ (this is the particular case of the comma
category): objects of $\left( {\A  \downarrow \C }
\right)$ are morphisms of $\C $ with source $\A $,
and morphisms of $\left( {\A  \downarrow \C }
\right)$ from one such object $F\colon \A  \to \B $ to
another $F'\colon \A  \to \mathbf{B'}$ are morphisms $G\colon
\B  \to \mathbf{B'}$ of $\C $ such that $F'=G \circ
F$.

Suppose $\C $ is a subcategory of $\PDE$, $\A $ is an object of $\C$. 
Then the category $\left( {\A  \downarrow \C } \right)$ of objects 
under $\A $ describes collection of quotient equations for $\A $ and 
their interconnection in the framework of $\C $.

To each morphism $F\colon \A  \to \B $ with source
$\A $ (that is to each object of the comma category $\left(
{\A  \downarrow \C } \right)$) assign the set
$F^{*}(\Sol (\B )) \subseteq \Sol (\A)$ of such solutions of $\A $ 
that ``projected'' onto underying space  of $\B$
(space of dependent and independent variables). 
We can identify such morphisms that generated the same sets of 
solutions of $\A$, that is identify isomorphic objects of the comma
category $\left( {\A  \downarrow \C } \right)$.

Describe the situation more explicitly. An equivalence class of
epimorphisms with source $\A $ is called a quotient object of
$\A $, where two epimorphisms $F\colon \A  \to
\B $ and ${F}'\colon\A  \to \B '$ are equivalent if
$F'=I \circ F$ for some isomorphism $I\colon \B  \to
\B '$ \cite{Maclain}.
If $F\colon \A  \to \B $ and ${F}'\colon\A  \to
\B '$ are equivalent, then they lead to the same subsets of
the solutions of $\A$: $F^{*}(\Sol(\B )) =
F'^{*}(\Sol(\mathbf{B'}))$. So if we interested only in the sets of
the solutions of $\A$, then all representatives of the
same quotient object have the same rights.

Therefore, the following problems naturally arise:
\begin{itemize}
  \item
  to study all morphisms with given source,
  \item
  to choose a "simplest" representative from every equivalence
class, or to choose representative with the simplest target (that is
the simplest quotient equation).
\end{itemize}
In order to deal with these problems, we develop a special-purpose language.

Let us introduce a number of partial orders on the class of all
categories to describe arising situations. First of all, we define a
few types of subcategories.

\begin{defi}\label{def_subcat}
Suppose $\C $ is a category, $\C _1$ is a subcategory of $\C $.
\begin{itemize}
\item
$\C _1$ is called a \textbf{wide} subcategory of $\C $ if
all objects of $\C $ are objects of $\C _1$.
\item
$\C _1$ is called a \textbf{full} subcategory of
$\C $ if every morphism in $\C $ with source and
target from $\C _1$ is a morphism in $\C _1$.
%New ===========================================
\item
We say that $\C _1$ is \textbf{full under isomorphisms in} $\C $ 
if every isomorphism in $\C $ with source and
target from $\C _1$ is an isomorphism in $\C _1$.
%New ===========================================
\item
We say that $\C _1$ is \textbf{closed in} $\C $ if
every morphism in $\C $ with source from $\C _1$ is
a morphism in $\C _1$. (Note that every subcategory that is
closed in $\C $ is full in $\C $.)
\item
We say that $\C _1$ is \textbf{closed under isomorphisms in}
$\C $ if every isomorphism in $\C $ with source from
$\C _1$ is an isomorphism in $\C _1$.
\item
We say that $\C _1$ is \textbf{dense in} $\C $ if
every object of $\C $ is isomorphic in $\C $ to an
object of $\C _1$.
\item
We say that $\C _1$ is \textbf{plentiful in} $\C $
if for every morphism $F\colon \A  \to \B $ in
$\C $, $\A  \in \Ob _{\C _1}$, there exist an
isomorphism $I\colon \B  \to \mathbf{C}$ in $\C $
such that $I \circ F \in \Hom _{\C _1}$ (in other words, for
every quotient object of $\A $ in $\C $ there exist a
representative of this quotient object in $\C _1$). Such
morphism $I \circ F$ we call $\C _1$\textbf{-canonical} for
$F$.
\item
We say that $\C _1$ is \textbf{fully dense (fully plentiful)
in} $\C $ if $\C _1$ is a full subcategory of
$\C $ and $\C _1$ is dense (plentiful) in
$\C $.
\end{itemize}
\end{defi}

\noindent The first two parts of the Definition are standard notions of
category theory, but the notions of the other parts are introduced here
for the sake of description of the $\PDE$ structure. 
%The notion ``closed under isomorphisms'' is introduced here for symmetry; it is not used in the next sections.

\begin{rem}
Using the notion of ``the category of objects under $\A $'',
we can define the notions of closed subcategory and plentiful
subcategory by the following way:
\begin{itemize}
\item
$\C _1$ is closed in $\C $ if for each $\A
\in \Ob _{\C _1}$ the category $\left( {\A
\downarrow \C _1} \right)$ is wide in $\left( {\A
\downarrow \C } \right)$.
\item
$\C _1$ is plentiful in $\C $ if for each
$\A  \in \Ob _{\C _1}$ the category $\left(
{\A  \downarrow \C _1} \right)$ is dense in $\left(
{\A  \downarrow \C } \right)$.
\end{itemize}
\end{rem}

\begin{rem}
$\C _1$ is fully dense in $\C $ if and only if the embedding functor 
$\C_1 \to \C$ defines an equivalence of these categories.
\end{rem}

Choose some category $\U$, which is big enough to contain all needful for us categories as it's subcategories. 
For our purposes $\U=\PDE$ will be sufficient.

Define the category $\Ug$ whose objects are subcategories $\C $ of $\U$,
and a collection $\Hom_{\U_{\geq}}\left( \C _1, \C _2\right)$ of morphisms from $\C _1$ to $\C _2$ is a
one-element set if $\C _2$ is subcategory of $\C _1$ and empty otherwise, so
an arrow from $\C _1$ to $\C _2$ in $\Ug$ means that $\C _2$ is the subcategory of $\C _1$.
Let $\U _{=}$ be the discrete wide subcategory of $\Ug$, that is objects of $\U _{=}$ are all subcategories $\C $ of $\U$,
and the only morphisms are identities, so $\C _1$ and $\C _2$ are connected by arrow in $\U _{=}$ only if $\C _1 = \C _2$.

\begin{defi}\label{def_cat_cap}
Suppose $\C _1$, $\C _2$ are subcategories of $\U$. 
A subcategory of $\U$, whose objects are objects of $\C _1$ and $\C _2$ simultaneously, 
and whose morphisms are morphisms of $\C _1$ and $\C _2$ simultaneously, 
is called an \textbf{intersection} of $\C _1$ and $\C _2$ 
and is denoted by $\C _1 \cap \C _2$.
\end{defi}

\noindent
In other words, $\C _1 \cap \C _2$ is the fibered sum of 
$\C _1$ and $\C _2$ in $\Ug$.

\begin{lem}\label{lem_cap12}
Suppose $\C _1$ is closed in $\C $, and
$\C _2$ is (full/closed/dense/plentiful) subcategory of $\C $;
then $\C _1 \cap \C _2$ is closed in $\C _2$
and is (full/closed/dense/plentiful) subcategory of $\C _1$.
\end{lem}

%\begin{table}[tbph]
\begin{floatingtable}[r]{
\includegraphics[scale=0.9]{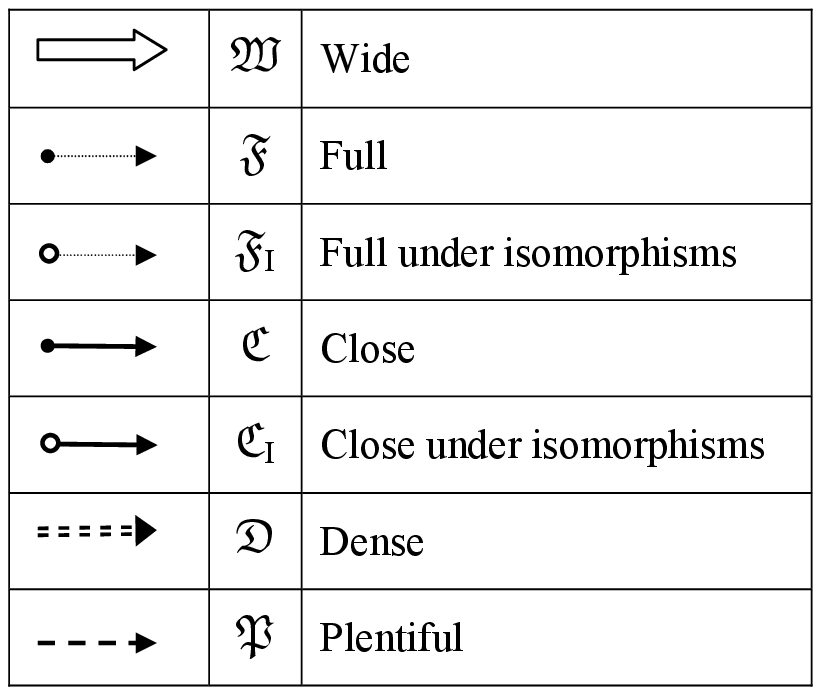}
}
\caption{Basic meta-categories (arrows)}
\label{TableArrows}
\end{floatingtable}
%\end{table}

Now we introduce some graphic designations for various types of
subcategories of $\Ug$ (see Table~\ref{TableArrows}). 
These designations will be used, particularly, for the representation 
of the structure of the category of parabolic equations described below.

We shall use the term ``meta-category'' both for the category
$\Ug$ and for its subcategories defined below to avoid
confusion between $\Ug$ and ``ordinary'' categories which
are objects of $\Ug$; and we shall use Gothic script for
meta-categories except $\Ug$. 
One may view these meta-categories as a partial
orders on the class of all subcategories of $\U$; 
we prefer category terminology here since this allow us to use category 
constructions for the interrelations between various partial orders.

Define $\mathfrak{W}$, $\mathfrak{F}$, $\mathfrak{F}_{\I}$, 
$\mathfrak{C}$, $\mathfrak{C}_{\I}$, $\mathfrak{D}$, and $\mathfrak{P}$ 
that are wide  subcategories of meta-category $Ug$. 
Objects of 
%$\mathfrak{W}$, $\mathfrak{F}$, $\mathfrak{F}_{\I}$, $\mathfrak{C}$, $\mathfrak{C}_{\I}$, $\mathfrak{D}$, and $\mathfrak{P}$ 
them are categories, but arrows from $\C _1$ to $\C _2$ have a different meaning:
\begin{itemize}
\item
in the meta-category $\mathfrak{W}$ it mean that $\C _2$
is wide subcategory of $\C _1$,
\item
in the meta-category $\mathfrak{F}$ it means that $\C _2$
is full subcategory of $\C _1$,
\item
in the meta-category $\mathfrak{F}_{\I}$ it means that $\C _2$
is full under isomorphisms in $\C _1$,
\item
in the meta-category $\mathfrak{C}$ it means that $\C _2$
is closed subcategory of $\C _1$,
\item
in the meta-category $\mathfrak{C}_{\I}$ it means that $\C _2$
is closed under isomorphisms in $\C _1$,
\item
in the meta-category $\mathfrak{D}$ it means that $\C _2$
is dense subcategory of $\C _1$,
\item
in the meta-category $\mathfrak{P}$ it means that $\C _2$
is plentiful subcategory of $\C _1$,
\end{itemize}
\noindent 
%Each of these types of meta-categories defines a certain partial order on the collection of all categories.
We shall denote the intersections of these meta-categories by the
concatenations of appropriate letters, for example:
$\mathfrak{FD}=\mathfrak{F} \cap \mathfrak{D}$.

\begin{lem}\label{lem_cap}
\noindent
\begin{itemize}
\item
$\mathfrak{F}_{\I} \cap \mathfrak{P} = \mathfrak{F} \cap \mathfrak{P}$,
\item
$\mathfrak{C}_{\I} \cap \mathfrak{P} = \mathfrak{C}$,
\item
$\mathfrak{F} \cap \mathfrak{P} \cap \mathfrak{D} = \mathfrak{F} \cap \mathfrak{D}$.
\end{itemize}
\end{lem}

\begin{figure}[tbh]
\begin{center}
\includegraphics[width=0.75\textwidth]{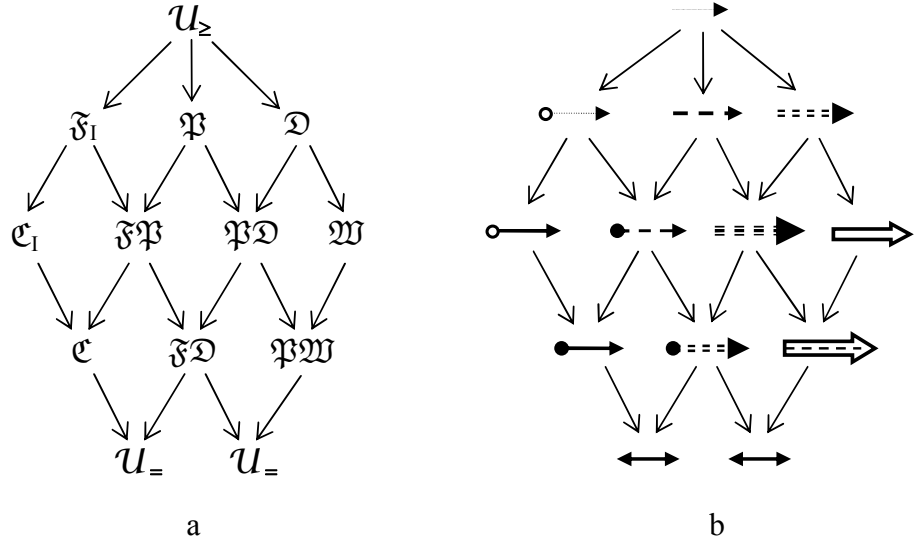}
\caption{Interrelations between basic meta-categories (arrows) and 
their intersections}
\label{FigArr}
\end{center}
\end{figure}

Interrelations between ``basic'' meta-categories 
$\mathfrak{W}$, $\mathfrak{F}$, $\mathfrak{F}_{\I}$, 
$\mathfrak{C}$, $\mathfrak{C}_{\I}$, $\mathfrak{D}$, $\mathfrak{P}$ 
and their intersections (``composed'' meta-categories) 
are represented on Fig.~\ref{FigArr}(a). 
Here an arrow means the predicate ``to be subcategory of''; 
we shall call it ``the meta-arrow''. 
For example, meta-arrow from $\mathfrak{D}$ to
$\mathfrak{W}$ means that $\mathfrak{W}$ is subcategory of
$\mathfrak{D}$. In the language of ``ordinary'' categories this
meta-arrow means that the statement ``$\C _2$ is wide in
$\C _1$'' implies that $\C _2$ is dense in
$\C _1$. Everywhere on Fig.~\ref{FigArr}(a) a pair of meta-arrows with
the same target means that this meta-category (target of these
meta-arrows) is the intersection of two ``top'' meta-categories
(sources of these meta-arrows). For example, $\mathfrak{FD} =
\mathfrak{FP} \cap \mathfrak{PD}$. 
%At the bottom segment of Fig.~\ref{FigArr}(a) the notation $\mathfrak{E}$ is used for the 
%discrete wide subcategory of $\mathfrak{S}$ (it is defined by the following way: 
%objects of this meta-category are all categories, and morphisms are identities in $\mathfrak{S}$, 
%so $\C _1$ and $\C _2$ are connected by arrow in this meta-category only if $\C _1 = \C _2$).

On Fig.~\ref{FigArr}(b) the same scheme is represented as on 
Fig.~\ref{FigArr}(a), but letter names are replaced by the arrows of 
various types.

Instead of the investigation of all or the simplest morphisms with
the given source, we want to introduce a certain structure in $\PDE$, 
so that the position of an object in it gives information
on the morphisms from the object and on the kind of the simplest
representatives of equivalence classes of the morphisms. We
construct this structure by means of choosing some subcategories in
$\PDE$ and connecting them by the arrows from Fig.~\ref{FigArr}(b). 
The part of obtained structure of the category of parabolic equations 
is shown on Fig.~\ref{FigBig}. 
We use this structure also in Section \ref{sec:reac-diff} to
investigate nonlinear reaction-diffusion equation.

\section{The category of parabolic equations}
\label{sec:PE}

Let us consider the class $P\left( {X,T,\Omega} \right)$ of
differential operators on a connected smooth manifold $X$, dependent
on the time $t$ as on parameter, that are of the form
\begin{gather}
Lu = \sum_{i,j} b^{ij}(t,x,u)u_{ij}
 + \sum_{i,j} c^{ij}(t,x,u)u_{i} u_{j}
 + \sum_{i} b^{i}(t,x,u)u_{i} + q(t,x,u), \notag
\\ \qquad
x \in X, \, t \in T, \, u \in \Omega \notag
\end{gather}
\noindent in some neighborhood of each point, in some (and then
arbitrary) local coordinates $\left( {x^{i}} \right)$ on $X$. Here
subscript $i$ denotes partial derivative with respect to $x^{i}$,
quadratic form $b^{ij} = b^{ji}$ is positive definite, $c^{ij} =
c^{ji}$. Both $T$ and $\Omega $ may be bounded, semibounded or
unbounded open intervals of $\R$.

\begin{defi}\label{def_PE}
\textbf{The category of parabolic equations of the second order}
$\PE$ is a subcategory of $\PDE$, whose objects are pairs 
$\A  = \left( {N,E} \right)$, $N = T \times X \times \Omega $, 
where $X$ is a connected smooth manifold, $T$ and $\Omega $ are open 
intervals, $E$ is an equation of the form $u_{t}= Lu$,
$L \in P\left( {X,T,\Omega} \right)$ (more exactly, $E$ is the 
extended version of an equation $u_{t}= Lu$, that is closed submanifold 
of $J^2_{n+1}(T \times X \times \Omega)$, $n=\dim X$).
\end{defi}

\begin{exm}\label{ex_sph_sol}
Let $\Phi _{k} (x)$, $x \in \R ^{3}- \{0 \}$ be a spherical harmonic
of the $k$-th order. Then the map $(t,x,u) \mapsto \left( {t,\left|
{x} \right|,{{u} \left/ {{\Phi _{k} (x)}} \right.}} \right)$ defined
the morphism in the category $\PE$ from the object
$\A $ corresponding to equation $u_{t} = \Delta u$ and $X =
\R ^{3} - \left\{ {0} \right\}$, $T = \Omega = \R$, to the object
$\mathbf{A'}$ corresponding to equation ${u}'_{{t}'} =
{u}'_{{x'}{x'}} - k\left( {k + 1} \right){x'}^{ - 2}{u}'$ and ${X'}
= \R _{+}$, ${T}' = {\Omega} ' = \R$. One may assign to the set
$\Sol(\mathbf{A'})$ of all solutions of the quotient equation the
set $F^{*}(\Sol(\mathbf{A'}))$ of such solutions of the original
equation that may be written in the form $u = \Phi _{k}
(x){u}'\left( {t,\left| {x} \right|} \right)$.
\end{exm}

\begin{exm}\label{ex_mor_to_me}
The following example shows that not every endomorphism in
$\PE$ is an automorphism. Consider object $\A $, for
which $X = S^{1} = \R \bmod 1$, $T = \Omega = \R $, $E\colon u_{t} =
u_{xx} $. Then morphism from $\A $ to $\A $ defined by
the map $(t,x,u) \mapsto \left( {4t,2x,u} \right)$ has no inverse.
\end{exm}

\begin{thm}\label{thm_mor_PE}
Every morphism in $\PE$ is of the form
\begin{equation}
\label{eqMorPE} (t,x,u) \mapsto \left( {{t}'(t),{x'}\left( {t,x}
\right),{u}'(t,x,u)} \right),
\end{equation}
\noindent with submersive ${t}'(t)$, ${x'}(t,x)$, ${u}'\left(
{t,x,u} \right)$. Isomorphisms in the category $\PE$ are exactly
diffeomorphisms of the form \eqref{eqMorPE}.
\end{thm}

\section{Certain subcategories of $\PE$: classification of the 
parabolic equations}
\label{sec:subcat}

\begin{figure}[tbph]
\begin{center}
\includegraphics[width=\textwidth]{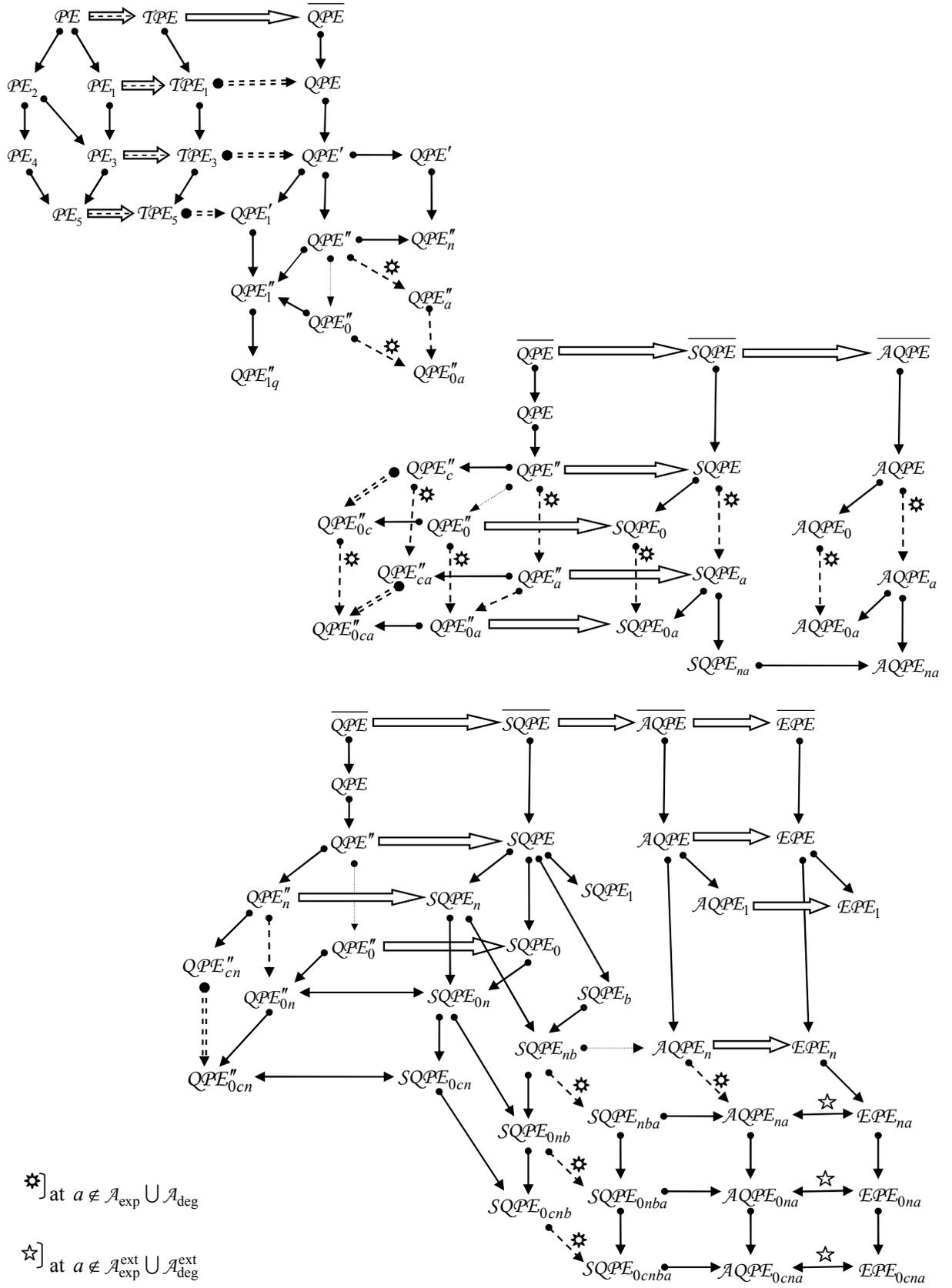}
\caption{The part of the structure of the category of parabolic equations}
\label{FigBig}
\end{center}
\end{figure}

Certain parts of the $\PE$ structure described below
are depicted schematically on Fig.~\ref{FigBig}.
The full structure is not depicted here in view of its awkwardness.

%\clearpage

Let us consider five full subcategories $\PE_{k}$ of the
category $\PE$, $1 \leq k \leq 5$, whose objects are
equations that can be written locally in the following form:
\begin{align}
 & u_{t} = \sum_{i,j} b^{ij}(t,x,u) 
   \left( {u_{ij} + \lambda (t,x,u)u_{i} u_{j}} \right)
   + \sum_{i} b^{i}(t,x,u)u_{i} + q(t,x,u)
\tag{$\PE_{1}$}\\
 & u_{t} = a(t,x,u) \sum_{i,j} \bar {b}^{ij} (t,x)u_{ij}
   + \sum_{i,j} c^{ij}(t,x,u)u_{i} u_{j}
   + \sum_{i} b^{i}(t,x,u)u_{i} + q(t,x,u)
\tag{$\PE_{2}$}\\
 & u_{t} = a(t,x,u) \sum_{i,j} \bar {b}^{ij}(t,x)
  \left( {u_{ij} + \lambda (t,x,u)u_{i} u_{j}} \right)
  + \sum_{i} b^{i} (t,x,u)u_{i} + q(t,x,u)
\tag{$\PE_{3}$}\\
 & u_{t} = \sum_{i,j} b^{ij}(t,x)u_{ij}
   + \sum_{i,j} c^{ij}(t,x,u)u_{i}u_{j}
   + \sum_{i} b^{i}(t,x,u)u_{i} + q(t,x,u)
\tag{$\PE_{4}$}\\
 & u_{t} = \sum_{i,j} b^{ij}(t,x)
   \left( {u_{ij} + \lambda (t,x,u)u_{i} u_{j}} \right)
   + \sum_{i} b^{i}(t,x,u)u_{i} + q (t,x,u)
\tag{$\PE_{5}$}
\end{align}

\begin{rem}
Everywhere in the paper we use notation of a category equipped with
a subscript and/or primes for its full subcategory. For example,
defined below $\QPE_{k}$, $\QPE'$ and
$\QPE'_{k}$ are full subcategories of $\QPE$.
\end{rem}

\begin{rem}\label{rem:a}
Note that in equations of the categories $\PE_2$ and $\PE_3$
a function $a$ is determined up to multiplication by arbitrary function
from $T \times X$ to $\R ^+$; moreover, it is determined only locally.
Nevertheless we can lead these equations to the equations of the same 
form but with globally defined function
$a\colon T \times X \times \Omega \to \R ^{+}$.
For example, we can require that conditions $a(t,x,u_0) \equiv 1$
take place, where $u_0$ is arbitrary point of $\Omega$.
Everywhere below we shall assume that function $a$ is globally determined
on $ T \times X \times \Omega$.
\end{rem}

\begin{thm}\label{thm_PE}
\noindent
\begin{enumerate}
\item
$\PE_{1} $ and $\PE_{2} $ are closed in $\PE$.
\item
$\PE_{3} =\PE_{1} \cap\PE_{2} $
is closed in $\PE_{1} $, in $\PE_{2}$, and in $\PE$.
\item
$\PE_{4} $ is closed in $\PE_{2} $ and in $\PE$.
\item
$\PE_{5} =\PE_{3} \cap\PE_{4} $
is closed in $\PE_{3} $, in $\PE_{4} $, and in $\PE$.
\end{enumerate}
\end{thm}

\begin{defi}\label{def_cat_overline}
Consider wide subcategories $\TPE$, $\overline {\QPE}$, 
$\overline {\SQPE} $, $\overline {\AQPE} $, and $\overline {\EPE} $ 
of $\PE$, whose morphisms are of the following form:
\[
(t,x,u) \to
\begin{cases}
 \left( {t,y(t,x),v\left( {t,x,u} \right)} \right)
& \text{for }\TPE \\
 \left( {t,y(t,x),\varphi (t,x)u + \psi (t,x)} \right)
& \text{for }\overline {\QPE} \\
 \left( {t,y(x),\varphi (t,x)u + \psi (t,x)} \right)
& \text{for }\overline {\SQPE} \\
 \left( {t,y(x),\varphi (x)u + \psi (x)} \right)
& \text{for }\overline {\AQPE} \\
 \left( {t,y(x),u} \right)
& \text{for }\overline {\EPE}
\end{cases}
\]
Denote $\TPE _k=\TPE \cap \PE_k$.
\end{defi}

\begin{thm}\label{thm_TPE}
\noindent
\begin{enumerate}
\item
$\TPE$ is wide and plentiful in $\PE$.
\item
$\TPE _k$ is closed in $\TPE$; it is wide and plentiful in $\PE_k$, 
$k=1..5$.
\end{enumerate}
\end{thm}

\begin{defi}\label{def_QPE}
Define the \textbf{category of quasilinear parabolic equations}
$\QPE$. $\QPE$ is the following full subcategory
of $\overline {{\QPE}} $: objects of $\QPE$ are
equations of the form
\[
u_{t} = \sum_{i,j} b^{ij}(t,x,u)u_{ij}
 + \sum_{i} b^{i}\left( {t,x,u} \right)u_{i} + q(t,x,u),
\tag{$\QPE$}\\
\]
\noindent
(in a local coordinates);
morphisms of $\QPE$ are maps of the form
\[
(t,x,u) \to  \left( {t,y(t,x),\varphi (t,x)u + \psi (t,x)} \right).
\]
\end{defi}

Denote by $\mathcal {A}_{nc} \left( {M,\Omega} \right)$ the set of
continuous positive functions $a\colon M \times \Omega  \to \R $ that
satisfies the condition
\[
\forall m \in M \quad \exists u_{1} ,u_{2} \quad
a\left( {m,u_{1}}  \right) \ne a\left( {m,u_{2}}  \right).
\tag{$\mathcal {A}_{nc}$}\\
\]
Define full subcategories of $\QPE$, whose objects are
equations of the following form:
\begin{align}
 & u_{t} = a(t,x,u) \sum_{i,j} \bar {b}^{ij}(t,x)u_{ij}
  + \sum_{i} b^{i}(t,x,u)u_{i} + q(t,x,u)
 \tag{$\QPE'$} \\
 & u_{t} = a(t,x,u) \sum_{i,j} \bar {b}^{ij}(t,x)u_{ij}
  + \sum_{i} b^{i}(t,x,u)u_{i} + q(t,x,u),
 \, a \in \mathcal {A}_{nc} \left( {T \times X} \right)
\tag{$\QPE'_n$} \\
 & u_{t} = \sum_{i,j} b^{ij}(t,x)u_{ij}
  + \sum_{i} b^{i}(t,x,u)u_{i} + q(t,x,u)
\tag{$\QPE'_1$} \\
 & u_{t} = a(t,x,u)\left( \sum_{i,j} \bar {b}^{ij}(t,x)u_{ij}
     + \sum_{i} \bar {b}^{i}(t,x)u_{i} \right) 
     + \sum_{i} \xi ^{i}(t,x)u_{i} + q(t,x,u)
\tag{$\QPE''$} \\
 & u_{t} = a(t,x,u )\left( \sum_{i,j} \bar {b}^{ij} (t,x)u_{ij}
   + \sum_{i} \bar {b}^{i}(t,x)u_{i} \right) + q(t,x,u)
\tag{$\QPE''_0$} \\
 & u_{t} = a(u) \left( \sum_{i,j} \bar {b}^{ij}(t,x)u_{ij}
   + \sum_{i} \bar {b}^{i}(t,x)u_{i} \right)
   + \sum_{i} \xi ^{i} (t,x)u_{i} + q(t,x,u)
\tag{$\QPE''_a(a)$} \\
 & u_{t} = \sum_{i,j} b^{ij}(t,x)u_{ij} + \sum_{i} b^{i}(t,x)u_{i}
  + q(t,x,u),
\tag{$\QPE''_1$} \\
 & u_{t} = \sum_{i,j} b^{ij}(t,x)u_{ij} + \sum_{i} b^{i}(t,x)u_{i}
  + q_{1} (t,x)u + q_{0} (t,x),
\tag{$\QPE''_{1q}$}
\end{align}
\noindent where $a>0$. The family of categories $\QPE''_{a}(a)$ is 
parameterized by functions $a\left( { \cdot}  \right)$, that is one 
assigns the category $\QPE''_{a} (a)$ to each continuous positive 
function $a \colon \Omega \to \R$.

Furthermore, consider the full subcategory $\QPE_{c}$ of
$\QPE$, whose objects are equations from $\QPE$
posed on a {\it compact} manifolds $X$.

Let us introduce the following notation for the intersections of 
enumerated ``basic'' subcategories:
for a string $\sigma$ we set $\QPE_{\sigma}  = \cap
\left\{ {\QPE_{\alpha}  \colon \alpha \in \sigma}  \right\}$,
$\QPE_{\sigma} ^{\beta}  =\QPE_{\sigma}  \cap
\QPE^{\beta} $. Particularly, $\QPE''_{0n}$
denotes the intersection $\QPE'_{n} \cap
\QPE''_{0}$.

In the same way as in Remark~\ref{rem:a}, if we impose the condition
$a(u_0)=1$, then we obtain the global function $a(u)$ for any
equation from $\QPE''_a(a)$; such function $a(u)$ does not
depend on the choice of neighborhood in $T \times X \times \Omega$ and 
on local coordinates.

\begin{thm}\label{thm_QPE}
\noindent
\begin{enumerate}
\item
$\QPE$ is closed in $\overline {{\QPE}} $
and is fully dense in $\TPE _1$.
\item
$\QPE_{c} $ is closed in $\QPE$.
\item
$\QPE' =\QPE \cap\PE_{2} =\QPE\cap\PE_{3} $
is fully dense in $\TPE _{3} $ and is closed in $\QPE$.
\item
$\QPE'_1 = \overline {{\QPE}} \cap\PE_{5} = \QPE'\cap\PE_{5}$
is fully dense in $\TPE _{5} $ and is closed in $\QPE'$.
\item
$\QPE''$ is closed in $\QPE'$.
\item
$\QPE''_{1} = \QPE'' \cap\PE_{5} =
\QPE'' \cap \QPE'_{1} = \QPE''_{a} \left( {1} \right)$
is closed in $\QPE'_{1} $, in $\QPE''$, and in $\QPE''_{0} $.
\item
$\QPE''_{1q} $ is closed in $\QPE''_{1} $.
\item
$\QPE'_{n} $ is closed in $\QPE'$.
\item
$\QPE''_{0n} $ is fully plentiful in $\QPE''_n$.
\item
$\QPE''_{0c} $ is fully dense in $\QPE''_c$.
\end{enumerate}
\end{thm}

Let $\mathcal{A}_{\exp} $ be the set of functions of the form
$a(u)=e^{\lambda u}H(u)$; let $\mathcal{A}_{\mathrm{deg}} $ be the
set of functions of the form $a(u) = \left( {u - u_{0}}
\right)^{\lambda} H\left( {\ln \left( {u - u_{0}}  \right)}
\right)$, where $\lambda $, $u_{0} $ are arbitrary constants and
$H\left( { \cdot}  \right)$ is arbitrary nonconstant periodic
function.

\begin{thm}\label{thm_QPE_a}
\noindent
\begin{enumerate}
\item
If $a \notin \mathcal{A}_{\exp} \cup \mathcal{A}_{\mathrm{deg}} $,
then $\QPE''_{a} (a)$ is fully plentiful in $\QPE''$.
\item
$\QPE''_{0a} (a)$ is fully plentiful in $\QPE''_{a} (a)$;
if $a \notin \mathcal{A}_{\exp} \cup \mathcal{A}_{\mathrm{deg}}$,
then $\QPE''_{0a} (a)$ is fully plentiful in $\QPE''_{0}$.
\item
$\QPE''_{0ca} (a)$ is fully dense in $\QPE''_{ca} (a)$.
\item
Suppose $\A $ is an object of $\QPE''_{a} (a)$,
$F\colon \A  \to \B $ is a morphism in $\PE$
such that there is no object of $\QPE''_{a} (a)$ isomorphic
to $\B $ in $\PE$ (that is $a(\cdot) \in
\mathcal{A}_{\exp} \cup \mathcal{A}_{\mathrm{deg}} $). Then there
exist an object of $\QPE''$ isomorphic to $\B $ such
that the composition of $F\colon \A  \to \B $ with
this isomorphism is of the form
\[
(t,x,u) \to
\begin{cases}
{\left( {t,y(t,x),u + \psi(t,x)} \right)},
 & {a \in \mathcal{A}_{\exp}}  \\
{\left( {t,y(t,x),v_0 + \left({u-u_0} \right)\exp\left( {\psi(t,x)} 
 \right)}\right)},
 & {a \in \mathcal{A}_{\mathrm{deg}}}  \\
\end{cases}
\]
\noindent In addition, for each $t\in T, x_1,x_2 \in X$ with
$y(t,x_1)=y(t,x_2)$ the difference $\psi(t,x_2)-\psi(t,x_1)$ is an
integral multiple of $\hat{H}$, where $\hat H$ is the period of
function $H$. The same assertion holds if we replace
$\QPE''_{a} (a)$ by $\QPE''_{0a} (a)$ and replace
$\QPE''$ by $\QPE''_{0} $.
\end{enumerate}
\end{thm}

\begin{exm}\label{ex_sin}
Consider the equation
\[
E\colon u_{t} = \left( {2 + \sin u} \right)u_{xx}
\]
\noindent
of the category $\QPE''_{0a}(f)$, where 
$X = T = \Omega = \R $, $f(u) = {2 + \sin u}$, $f \in \mathcal{A}_{\exp}$.
It admits both maps $(t,x,u) \mapsto
\left( {t,x\bmod {2\pi} ,u} \right)$ and $(t,x,u) \mapsto \left(
{t,x\bmod {2\pi},u + x} \right)$. In both cases $Y = S^{1}$.
In the first case the quotient equation is of the form $v_{t} =
 \left( {2 +\sin v} \right)v_{yy} $
and is an object of $\QPE''_{0a} (f)$.
In the second case the quotient equation is of the form 
$v_{t} = \left( {2 + \sin(v+y)} \right)v_{yy}$;
it is the object of $\QPE''_{0}$, but is not isomorphic to any object 
of $\QPE''_{0a} (f)$.
\end{exm}

\begin{defi}\label{def_SQPE}
The \textbf{category of semi-autonomous quasilinear parabolic equations} 
$\SQPE$ is the intersection $\overline {{\SQPE}} \cap \QPE''$.
In other words, $\SQPE$ is the full subcategory of 
$\overline {{\SQPE}}$ and the wide subcategory of $\QPE''$, 
whose objects are equations of the form
\[
u_{t} = a(t,x,u)\left( \sum_{i,j} \bar {b}^{ij} (t,x)u_{ij}
 + \sum_{i} \bar {b}^{i}(t,x)u_{i} \right) + \sum_{i} \xi ^{i}(t,x)u_{i}
 + q(t,x,u),
\tag{$\SQPE$}
\]
\noindent and morphisms are maps of the form 
$(t,x,u) \mapsto \left( {t,y(x),\varphi (t,x)u + \psi (t,x)} \right)$.

Define the following full subcategories of $\SQPE$: 
$\SQPE_{\sigma} = \SQPE \cap \QPE''_{\sigma}$; $\SQPE_{b}$ 
is the category, whose objects are equations of the form
\[
 u_{t} = a(t,x,u)\left( \sum_{i,j} {\bar {b}^{ij}(x)u_{ij}
 + \sum_{i} \bar {b}^{i}(t,x)u_{i}}  \right) 
 + \sum_{i} \xi ^{i} (t,x)u_{i} + q(t,x,u).
\tag{$\SQPE_b$}
\]
\end{defi}

\begin{thm}\label{thm_SQPE}
\noindent
\begin{enumerate}
\item
$\SQPE$ is closed in $\overline {{\SQPE}}$.
\item
$\SQPE_{0} = \overline {{\SQPE}} \cap \QPE''_{0} $, $\SQPE_{n} = 
\overline {{\SQPE}} \cap \QPE''_{n} $, and $\SQPE_{b} $
are closed in  $\SQPE$.
\item
$\SQPE_{0n} $ coincides with $\QPE''_{0n} $; it is
closed in $\SQPE_{0}$ and in $\SQPE_{n}$.
\item
$\SQPE_{1} = \overline {{\SQPE}} \cap \QPE''_{1} =
\SQPE_{a} \left( {1} \right)$
is closed in $\SQPE_{0} $.
\item
If $a \notin \mathcal{A}_{\exp} \cup \mathcal{A}_{\mathrm{deg}}$,
then $\SQPE_{a} (a)$ is fully plentiful in $\SQPE$.
\end{enumerate}
\end{thm}

\begin{defi}\label{def_AQPE}
The \textbf{category of autonomous quasilinear parabolic equations}
$\AQPE$ is the full subcategory of $\overline {{\AQPE}} $,
each object of $\AQPE$ is an equation of the form
\[
u_{t} = a(x,u)\left( {\Delta u + \eta \nabla u} \right) + \xi \nabla u
 + q(x,u)
\tag{$\AQPE$}
\]
\noindent
posed on a Riemann manifold $X$ equipped with vector fields $\xi$, $\eta$.

\noindent Define the following full subcategories of
$\AQPE$: $\AQPE_{\sigma} = \AQPE \cap
\QPE''_{\sigma}$ is the category, whose objects are
equations of the form
\begin{align}
 & u_{t} = a(x,u)\left( {\Delta u + \eta \nabla u} \right) +
  \xi \nabla u + q(x,u),
 \quad a \in \mathcal {A}_{nc} (X),
\tag{$\AQPE_n$} \\
 & u_{t} = a(x,u) (\Delta u + \eta \nabla u) + q(x,u),
\tag{$\AQPE_0$} \\
 & u_{t} = a(u) (\Delta u + \eta \nabla u) + \xi \nabla u + q(x,u),
\tag{$\AQPE_a(a)$} \\
 & u_{t} = \Delta u + \xi \nabla u + q(x,u).
\tag{$\AQPE_1$}
\end{align}
\end{defi}

\begin{thm}\label{thm_AQPE}
\noindent
\begin{enumerate}
\item
$\AQPE$ is closed in $\overline {{\AQPE}} $.
\item
$\AQPE_{n} $ is closed in $\AQPE$ and full in $\SQPE_{bn}$.
\item
$\AQPE_{0} $ and $\AQPE_{1} $ are closed in $\AQPE$.
\item
If $a\left( { \cdot}  \right) \notin \mathcal{A}_{\exp}
 \cup \mathcal{A}_{\mathrm{deg}}$,
then $\AQPE_{a} (a)$ is fully plentiful in $\AQPE$.
\item
$\AQPE_{na} (a)$ is closed in $\SQPE_{na} (a)$.
\end{enumerate}
\end{thm}

\begin{defi}\label{def_EPE}
Define the following full subcategories of the category
$\overline {{\EPE}}$ (its morphisms are a maps of the form
$(t,x,u) \mapsto (t,y(x),u)$):
\begin{gather*}
\EPE = \overline {{\EPE}} \cap \AQPE,  \notag \\
\EPE_{\sigma}  = \overline {{\EPE}} \cap \AQPE_{\sigma}, \notag \\
\EPE_{a} (a) = \overline {{\EPE}} \cap \AQPE_{a} (a).
\end{gather*}
\end{defi}

\begin{thm}\label{thm_EPE}
\noindent
\begin{enumerate}
\item
$\EPE$ is closed in $\overline {\EPE}$ and wide in
$\AQPE$.
\item
$\EPE_{n} $, $\EPE_{0} $, $\EPE_{1} $, and
$\EPE_{a} (a)$
are closed in $\EPE$.
\item
If
$a \notin \mathcal{A}_{\exp}^{\ext}
 \cup\mathcal{A}_{\mathrm{deg}}^{\ext}$,
then $\EPE_{a} (a)$ coincides with $\AQPE_{a} (a)$.
\end{enumerate}
\end{thm}

\noindent Here $\mathcal{A}_{\exp}^{\ext} $ is the set of
functions $a(u)$ of the form $a(u) = e^{\lambda u}H(u)$,
$\mathcal{A}_{\mathrm{deg}}^{\ext} $ is the set of functions
of the form $a(u) = \left({u-u_0} \right)^{\lambda} H\left( {\ln
\left({u-u_0} \right)} \right)$, $\lambda $, $u_{0} $ are arbitrary
constants, $H (\cdot)$ is arbitrary periodic function (that is
$\mathcal{A}_{\exp} \subset \mathcal{A}_{\exp}^{\ext}$,
$\mathcal{A}_{\deg} \subset \mathcal{A}_{\deg}^{\ext}$).

\begin{figure}[tbh]
\begin{center}
\includegraphics
%[scale=0.6]
[width=0.9\textwidth]
{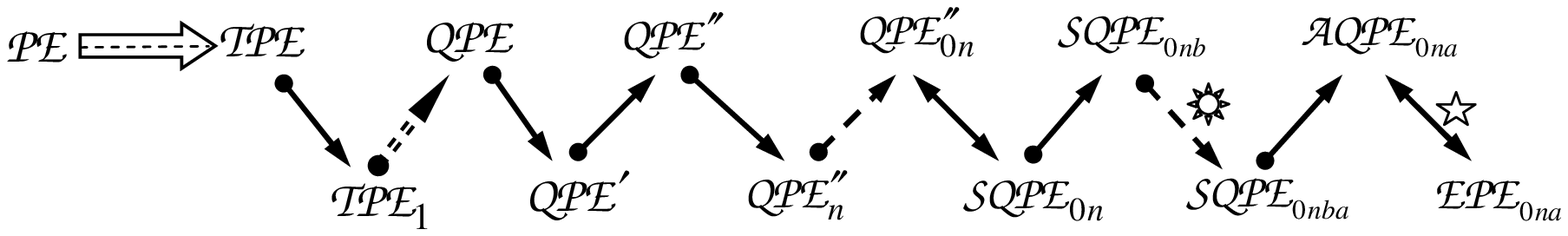}
\caption{The sequence of arrows from $\PE$ to $\EPE_{0na} (a)$}
\label{FigSeq}
\end{center}
\end{figure}

Let us consider the sequence depicted on Fig.~\ref{FigSeq}. 
Selecting the ``weakest'' arrow in this sequence, we obtain 
the following result.

\begin{cor}\label{cor:seq}
\noindent
\begin{enumerate}

\item
If $a \notin \mathcal{A}_{\exp} \cup \mathcal{A}_{\mathrm{deg}} $,
$a \ne \mathrm{const}$, then $\AQPE_{0a} (a)$ is fully
plentiful in $\TPE $ and plentiful in $\PE$.
\item
If $a \notin \mathcal{A}_{\exp}^{\ext}\cup\mathcal{A}
_{\mathrm{deg}}^{\ext}$,
$a \ne \mathrm{const}$, then $\EPE_{0a} (a)$ is fully
plentiful in $\TPE $ and plentiful in $\PE$.
\end{enumerate}
\end{cor}

\section{Factorization of the reaction-diffusion equation}
\label{sec:reac-diff}
Let us consider a nonlinear reaction-diffusion equation
\[
u_{t} = a(u)\left( {\Delta u + \eta \nabla u} \right) + q(x,u),
\]
\noindent $a \ne \mathrm{const}$, 
posed on a Riemann manifold $X$ equipped with a vector field $\xi$.
This equation defines the object $\A $ of the category
$\EPE_{0na} (a)$. Using Corollary~\ref{cor:seq}, we get the following
result:

\begin{cor}
\noindent
\begin{enumerate}
\item
If $a \notin \mathcal{A}_{\exp} \cup \mathcal{A}_{\mathrm{deg}} $,
then for every morphism $F\colon \A  \to \B $
of the category $\PE$
there exist isomorphism $I\colon \B  \to \B '$
of the form \eqref{eqMorPE}
(in other words, bijective change of variables in the quotient equation),
transforming $F$ to the canonical morphism of the form
\begin{equation}
\label{eqCanMor}
 (t,x,u) \to
\begin{cases}
  \left( {t,y(x),\varphi (x)u + \psi (x)} \right)
& \text{at } a \notin \mathcal{A}_{\exp} \cup \mathcal{A}_{\mathrm{deg}};
\\
  \left( {t,y(x),u} \right)
& \text{at } a \notin  \mathcal{A}_{\exp}^{\ext}
 \cup \mathcal{A}_{\mathrm{deg}}^{\ext}.
\end{cases}
\end{equation}
The corresponding canonical quotient equation ${B}'$ posed on the
Riemann manifold $Y$ is of the form
\begin{equation}
\label{eqEPE_end}
v_{t} = a(v)\left( {\Delta v + {\rm H}\nabla v} \right) +
 Q\left( {y,v} \right),
\quad y \in Y,
\quad {\rm H} \in TY
\end{equation}
\noindent
with the same function $a$.
\item
If a morphism $F\colon \A  \to \B $ of the category
$\TPE$ transforms $\A $ to an equation of the form
\eqref{eqEPE_end}, then $F$ is of the form \eqref{eqCanMor}.
\end{enumerate}
\end{cor}

\section{Comparison with the reduction by a symmetry group}
\label{sec:group}
As Remark~\ref{rem:FactByGroup} shows, our definition of morphism in $\PDE$ 
is a generalization of the reduction by a symmetry group. So we could
obtain the sets of solutions being more common than the sets of
group-invariant solutions of group analysis of PDE (though our
approach is more laborious owing to the non-linearity of the system of PDE 
describing a morphisms). In what follows we show this on the example of
a primitive morphism.

\begin{defi}\label{def_prim_mor}
A morphism $F\colon \A  \to \B $ of a category
$\C $ is called a \textbf{reducible in} $\C $ if in
$\C $ there are exists non-invertible morphisms $G\colon
\A  \to \mathbf{C}$, $H\colon \mathbf{C} \to \B $ such
that $F=H \circ G$. Otherwise, a morphism is called
\textbf{primitive in} $\C $.
\end{defi}

Note that the reduction of PDE by a symmetry group defines a
primitive morphism if and only if the group has no proper subgroups,
i.e. the group is a discrete cyclic group of prime order. The
reduction by any symmetry group that is not a discrete cyclic group 
of prime order may be always represented as a superposition of 
two nontrivial reductions, so the corresponding
morphism is a superposition of two non-invertible morphisms and
therefore is reducible. In particular, this situation takes place
for any nontrivial connected Lie group.

However, the situation for morphisms is completely different.
Even a morphism that decreases the number of independent variables
by 2 or more may be primitive; below we present an example of such a
morphism. However, in the Lie group analysis we always have
one-parameter subgroups of a symmetry group, so the morphism,
corresponding to a symmetry group, is always reducible.

\begin{exm}\label{ex_prim_mor}
Consider the following morphism $F\colon \A  \to \B $ in $\PE$:

\begin{itemize}

\item
$\A $ is heat equation $u_t=a(u)\Delta u$ posed on
$X = \left\{ {\left( {x,y,z,w} \right)\colon z < w} \right\} \subset
 \R ^{4}$ equipped with the metric
\[
g_{ij} =
\begin{pmatrix}
1 & 0 & 0 & 0 \\
0 & \gamma & \alpha & \beta \\
0 & \alpha & 1 & 0 \\
0 & \beta  & 0 & 1
\end{pmatrix},
\]
\noindent
where $\alpha = xe^{w}$, $\beta = xe^{z}$, $\gamma=1 +
\alpha ^{2} + \beta ^{2}$, $a \notin \mathcal{A}_{\exp} \cup
\mathcal{A}_{\mathrm{deg}}$. In the coordinate form $\A $
looks as
\begin{gather}
a^{ - 1}(u)u_{t} = u_{xx} + u_{yy} - 2\alpha u_{yz} - 2\beta u_{yw}
+ \left( {1 + \alpha ^{2}} \right)u_{zz} + \notag
\\ \qquad
 + 2\alpha \beta u_{zw} + \left( {1 + \beta ^{2}} \right)u_{ww} +
\left( {\alpha \beta}  \right)_{w} u_{z} + \left( {\alpha \beta}
\right)_{z} u_{w}. \notag
\end{gather}
\item
$\B $ is heat equation $a^{ - 1}(u)u_{t} = u_{xx} + u_{yy} $
posed on $Y = \left\{ {\left( {x,y} \right)} \right\} = \R ^{2}$
equipped with Euclidean metric.
\item
The morphism $F$ is defined by the map $\left( {t,(x,y,z,w),u}
\right) \mapsto \left( {t,(x,y),u} \right)$.
\end{itemize}

\noindent This morphism decreases the number of independent
variables by 2 and nevertheless is primitive in $\PE$.
\end{exm}

Additional examples of morphisms that are not defined by any
symmetry group of the given PDE, and also a detailed
investigation of the case $\dim Y = \dim X-1$, may be found in
\cite{Prokh2000}, \cite{Prokh2001}.

\section{Proofs}\label{sec:Proofs}

%\begin{proof}[Proof of Theorem \ref{thm_mor_PE}]
\subsection*{Proof of Theorem \ref{thm_mor_PE}}

Passing from the equation $u_{t} = Lu$ to the equation in the
extended jet bundle for unknown submanifold $L \subset X \times
T \times \Omega $ locally defined by the formula $f\left( {t,x,u} \right) =
0$, and expressing the derivatives of $u$ by the corresponding
derivatives of $f$, we obtain the following extended version of the
equation $E$:
\begin{equation}
\label{eq_f}
f_{t} f_{u}^{2} = \sum_{i,j} b^{ij}\left( {f_{ij} f_{u}^{2}
   - \left( {f_{iu} f_{j} + f_{ju} f_{i}}  \right)f_{u} 
   + f_{i} f_{j} f_{uu}} \right) 
   - \sum_{i,j} c^{ij}f_{i} f_{j} f_{u} 
   + \sum_{i} b^{i}f_{i} f_{u}^{2} - qf_{u}^{3}
\end{equation}
Suppose $F\colon \A  \to \A '$ is a morphism in $\PE$, 
${N}' = {X'} \times {T'} \times {\Omega} '$, 
${E}'$ is defined by the equation
\[
{u'} = \sum_{i',j'} {B}^{{i}'{j}'}\left( {{t'},{x'},{u'}} \right)
    {u'}_{{i}'{j}'} 
    + \sum_{i',j'} C^{{i}'{j}'}\left( {{t'},{x'},{u'}} \right)
    {u'}_{{i}'} {u'}_{{j}'} 
    + \sum_{i'} {B}^{{i}'}\left( {{t'},{x'},{u'}} \right){u'}_{{i}'} 
    + Q\left( {{t'},{x'},{u'}} \right).
\]
Consider the extended analog of the last equation:
\begin{multline}\label{eq_ff}
{f'}_{{t'}} {f'}_{{u'}}^{2} =
  \sum_{i',j'} B^{{i}'{j}'}\left( {{f'}_{{i}'{j}'} {f'}_{{u'}}^{2} 
  - \left( {{f'}_{{i}'{u'}} {f'}_{{j}'} 
  + {f'}_{{j}'{u'}} {f'}_{{i}'}} \right){f'}_{{u'}} 
  + {f'}_{{i}'} {f'}_{{j}'} {f'}_{{u'}{u'}}}
\right) - %\notag
\\ %\qquad
- \sum_{i',j'} C^{{i}'{j}'}{f'}_{{i}'} {f'}_{{j}'} {f'}_{{u'}} 
 + \sum_{i'} B^{{i}'}{f'}_{{i}'} {f'}_{{u'}}^{2} - Q{f'}_{{u'}}^{3} ,
\end{multline}
\noindent where the equation $f'\left( {t',x',u'} \right) = 0$
locally defines a submanifold $L '$ of $N'$.

Recall that $F\colon \left( {t,x,u} \right) \mapsto \left(
{{t'},{x'},{u'}} \right)$ is morphism in $\PE$ if and only
if for each point $\vartheta \in N$ and for each submanifold $L'$ of $N'$, 
$F(\vartheta) \in L '$, the following two
conditions are equivalent:

\begin{itemize}
\item the 2-jet of $L '$ at the point $F(\vartheta)$ satisfies
\eqref{eq_ff}
\item the 2-jet of $F^{-1}\left(L'\right)$ at the point $\vartheta$ satisfies 
\eqref{eq_f}.
\end{itemize}
\noindent In other words, conditions ``${f'}$ is solution of
\eqref{eq_ff}'' and ``$f$ is solution of \eqref{eq_f}'' must be
equivalent when
\[
f\left( {t,x,u} \right) = {f'}\left(
{{t'}\left( {t,x,u} \right),{x'}\left( {t,x,u} \right),{u'}\left(
{t,x,u} \right)} \right).
\]
\noindent
To find all such maps we use the
following procedure:

\begin{enumerate}

\item
Express derivatives of $f$ in \eqref{eq_f} through derivatives of ${f'}$:
\[
\frac{{\partial f}}{{\partial t}} =
\frac{{\partial {f'}}}{{\partial {t'}}}
\frac{{\partial {t'}}}{{\partial t}}
+ \frac{{\partial {f'}}}{{\partial {x'}^{{i}'}}}
\frac{{\partial {x'}^{{i}'}}}{{\partial t}}
+ \frac{{\partial {f'}}}{{\partial {u'}}}
\frac{{\partial {u'}}}{{\partial t}}
\]
\noindent
and so on.
\item
In the obtained identity substitute the combinations of the
derivatives of ${f'}$ for ${{\partial {f'}}}/{{\partial {t'}}}$ by
formula (\ref{eq_ff}). Then repeat this step for ${{\partial
^{2}{f'}}}/{{\partial {t'}^{2}}}$ in order to eliminate all
derivatives with respect to $t'$. After reducing to common
denominator, obtained identity will be of the form $\Phi =0$, where
$\Phi$ is rational function of partial derivatives of ${f'}$ with
respect to ${x'}$ and ${u'}$. The coefficients $\phi _{1} , \ldots
,\phi _{s}$ of $\Phi$ are functions of 4-jet of the map $F$.
\item
Solve the system of PDEs $\phi _{1} = 0, \ldots ,\phi _{s} = 0$ for
the map $F$.
\end{enumerate}

Let us realize this procedure. Note that we shall not write out
function $\Phi $ completely, but consider only some of its
coefficients. We shall use the obtained information about $F$ in
order to simplify $\Phi $ step by step in the following manner.

First note that derivatives of the forth order arise only in the term
${{\partial ^{2}{f'}}}/{{\partial {t'}^{2}}}$
when we fulfill step 2 of the above procedure.
Write this term before the final realization of step 2 for the sake of 
simplicity:
\begin{gather}
\Phi = \sum_{i,j} b^{ij}\left( {{t'}_{i} {t'}_{j} f_{u}^{2}
   - {t'}_{i} {t'}_{u} f_{j} f_{u} - {t'}_{j} {t'}_{u} f_{i} f_{u} 
   + {t'}_{u}^{2} f_{i} f_{j}} \right)
   \frac{{\partial ^{2}{f'}}}{{\partial {t'}^{2}}} + \ldots = 
\notag
\\ \qquad
= \sum_{i,j} b^{ij}\left( {{t'}_{i} f_{u} - {t'}_{u} f_{i}}  \right)
   \left( {{t'}_{j} f_{u} - {t'}_{u} f_{j}}  \right)
   \frac{{\partial ^{2}{f'}}}{{\partial {t'}^{2}}} + \ldots \notag
\end{gather}
\noindent The coefficient at ${{\partial ^{2}{f'}}}/{{\partial
{t'}^{2}}}$ must be zero, and the quadratic form $b^{ij}$ is
positive defined. We get ${t'}_{i} f_{u} = {t'}_{u} f_{i} $, so
\[{t'}_{i} \left( {{f'}_{{t'}} {t'}_{u} + {f'}_{{x'}} {x'}_{u} +
{f'}_{{u'}} {u'}_{u}}  \right) = {t'}_{u} \left( {{f'}_{{t'}}
{t'}_{i} + {f'}_{{x'}} {x'}_{i} + {f'}_{{u'}} {u'}_{i}}  \right)
\]
\noindent
(here and below we use notations
\[
{f'}_{{i}'} = \frac{\partial {f'}}{\partial {x'}_{{i}'}}, \quad
{f'}_{{x'}} {x'}_{u} =
\sum\limits_{{j}'} {{f'}_{{j}'} {x'}_{u}^{{j}'}}
\]
\noindent
and so on).
Hence we obtain the following system of equations:
\begin{equation}
\label{eq_thm1_tu}
\left\{
 {\begin{array}{l}
 {{t'}_{u} {u'}_{i} = {u'}_{u} {t'}_{i}}  \\
 {{t'}_{u} {x'}_{i} = {x'}_{u} {t'}_{i}}
 \end{array}}
\right.
\end{equation}

One of the following three conditions holds:

\begin{enumerate}

\item
${t'}_u = 0$, ${t'}_x = 0$;

\item
${t'}_u = 0$, ${t'}_x \ne 0$;
\item
${t'}_u \ne 0$.
\end{enumerate}

In the second case ${u'}_{u} = {x'}_{u} = 0$. Taking into account
equality ${t'}_u = 0$, we obtain a desired contradiction to the
assumption that $F$ is a submersion.

In the third case we get from (\ref{eq_thm1_tu}) 
${t'}_{x} = \omega {t'}_{u} $, 
${u'}_{x} = \omega {u'}_{u} $, 
${x'}^{{i}'}_{x} = \omega {x'}_{u}^{{i}'} $, 
where $\omega = {{t'}_{x}}/{{t'}_{u}} \in \Gamma\left(\pi^{*}T^{*} X\right)$, 
$\pi \colon N = T \times X \times \Omega \to X$ is the natural projection,
$\pi^{*}T^{*} X$ is the vector bundle over $N$ induced by $\pi$ from 
the cotangent bundle $T^{*} X$ over $X$,
$\omega = \sum_{i}\omega_i (t,x,u) \d x^i$ in local coordinates. 
This implies that $f_{x} = \omega f_{u} $. 
Substituting last formula to (\ref{eq_f}), we get
\[
f_{t} = f_{u} \left[ 
  \sum_{i,j} b^{ij} \left( 
     \frac{\partial\omega_i}{\partial x^j} 
     - \omega_j \frac{\partial\omega_i}{\partial u} 
    \right)
   - \sum_{i,j} c^{ij}\omega_i \omega_j 
   + \sum_{i} b^{i}\omega_i - q \right].
\]
Denote by $\zeta \left( {t,x,u} \right)$ the expression in square
brackets. Then $f_{t} = \zeta f_{u} $. Expressing derivatives of $f$
in terms of derivatives of ${f'}$, we get ${t'}_{t} = \zeta {t'}_{u}
$, ${x'}_{t} = \zeta {x'}_{u} $, ${u'}_{t} = \zeta {u'}_{u} $.
Consider the field of hyperplanes that kill the 1-form $\d t'$ in the
tangent bundle $TM$ (recall that ${t'}_u \neq 0$, so $\d t'$ is
nondegenerated). The differential of the map $F$ vanishes on these
hyperplanes, because $\d u' \wedge \d t' = \d x'^{i'} \wedge \d t' = 0$.
Therefore ${\rm rang} (\d F) \le 1$ and $F$ could not be submersive,
because $\dim N' \ge 3$.

Finally, only the first case is possible. Hence ${t'}$ is a function
of $t$, and ${f'}_{{t'}} $ could appear only in the representation
of $f_{t}$. Let us look at the terms of $\Phi $, containing 
$\left( {{f'}_{{u'}}} \right)^{ - 2}$:
\[
\Phi = \sum_{i',j',k',l'} {t'}_{t} {x'}_{u}^{{i}'} {x'}_{u}^{{j}'} 
{B'}^{{k}'{l}'}{f'}_{{i}'} {f'}_{{j}'} {f'}_{{k}'} {f'}_{{l}'}
{f'}_{{u'}{u'}} \left( {{f'}_{{u'}}} \right)^{ - 2} + \ldots .
\]
\noindent 
Substitution of any covector 
$\omega = \sum_{i'}\omega_{i'} \d {x'}^{i'} \in \Gamma\left(T^{*} X'\right)$
to the expression 
\[
\sum_{i',j',k',l'} {t'}_{t} {x'}_{u}^{{i}'} {x'}_{u}^{{j}'} 
{B'}^{k'l'}{\omega}_{i'} {\omega}_{j'} {\omega}_{k'} {\omega}_{l'}
= {t'}_{t} \left( \sum_{i'} {x'}_{u}^{i'}\omega_{i'} \right)^2
\left( \sum_{k',l'} {B'}^{{k}'{l}'}\omega_{k'} \omega_{l'} \right)
\]	
should give zero. 
The form ${B'}^{{k}'{l}'}$ is positively defined so 
$\sum_{k',l'} {B'}^{{k}'{l}'}\omega_{k'} \omega_{l'} >0$
when $\omega \neq 0$.
Taking into account that $F$ is submersive, we obtain ${t'}_{t} \ne 0$, 
so $\sum_{i'} {x'}_{u}^{i'}\omega_{i'} =0$ for any $\omega$, 
that is ${x'}_{u} \equiv 0$. 
Hence ${x'} = {x'}\left( {t,x} \right)$, ${t'} = {t'}\left( {t} \right)$.
\qed
%\end{proof}

\subsection*{Proof of Theorem \ref{thm_PE}}
%\pf[of Theorem \ref{thm_PE}]

%\subsection*{}

The map $\left( {t,x,u} \right) \mapsto \left( {\tau (t),y\left(
{t,x} \right),v(t,x,u)} \right)$ is a morphism in $\PE$ if
and only if
\begin{equation}
\label{eq_PE}
\left\{
\begin{aligned}
\tau _{t} B^{kl}
& = \sum_{i,j} b^{ij}y_{i}^{k} y_{j}^{l}  \\
\tau _{t} C^{kl}
&  =   \left( {\ln  U_{v}}  \right)_{v} B^{kl} 
  + U_{v} \sum_{i,j} c^{ij}y_{i}^{k} y_{j}^{l} \\
\tau _{t} B^{k}
& =  \sum_{i,j} b^{ij}y_{ij}^{k} + 2 \sum_{i,j} b^{ij} 
       \left( {\ln  U_{v}} \right)_{j} y_{i}^{k}
   + 2 \sum_{i,j} c^{ij}U_{j} y_{i}^{k} 
   + \sum_{i} b^{i}y_{i}^{k} - y_{t}^{k} \\
\tau _{t} Q
&  = U_{v}^{-1} \left( \sum_{i,j} b^{ij}U_{ij} 
    + \sum_{i,j} c^{ij}U_{i} U_{j}
    + \sum_{i} b^{i}U_{i} + q(t,x,U) - U_{t}  \right)
\end{aligned}
\right.
\end{equation}
\noindent where function $u = U\left( {t,x,v} \right)$ is the
inverse of the $v\left( {t,x,u} \right)$. 
The quotient equation is written as 
$v_{\tau} = \sum_{k,l} B^{kl}v_{kl} + \sum_{k,l} C^{kl}v_{k} v_{l} +
\sum_{k} B^{k}v_{k} + Q$. 
Here and below indexes $i$, $j$ relate to $x$, indexes $k$, $l$ relate to $y$.

By definition, all $\PE_{k} $ are full subcategories of $\PE$.

%\subsection*{}
\bigskip\noindent
1. Let us prove that $\PE_{1}$ is closed in $\PE$.
Suppose $\A  \in \Ob _{\PE_{1}} $, $F\colon \A  \to \B $ is a morphism 
in $\PE$. Then $c^{ij} = \lambda (t,x,u)b^{ij}$. 
From the second equation of system \eqref{eq_PE} we get
\[
C^{kl}\left( {\tau ,y,v} \right) = B^{kl}\left( {\tau ,y,v} \right)
\left[ {\tau _{t}^{ - 1} \left( {\ln  U_{v}}  \right)_{v} +
 \lambda \left( {t,x,u}\right)U_{v}}  \right].
\]
\noindent
The quadratic form $B^{kl}$ is nondegenerated at any point 
$\left( {\tau ,y,v}\right)$, so expression in square brackets is 
function of $\left( {\tau ,y,v} \right)$: 
$\tau _{t}^{ - 1} \left( {\ln  U_{v}} \right)_{v} 
 + \lambda (t,x,u)U_{v} = \Lambda \left( {\tau ,y,v}
\right)$, and $C^{kl}\left( {\tau ,y,v} \right) = \Lambda \left(
{\tau ,y,v} \right)B^{kl}(\tau ,y,v)$. Thus $\B  \in \Ob
_{\PE_{1}} $.

Let us show that $\PE_{2} $ is closed in $\PE$.
Suppose $\A  \in \Ob _{\PE_{2}}$, $F\colon
\A  \to \B $ is a morphism in $\PE$. Then
$b^{ij} = a(t,x,u)\bar {b}^{ij}(t,x)$. Using the first equation of
the system \eqref{eq_PE}, we receive
\[
\tau _{t} B^{kl} =
a(t,x,u)\left( \sum_{i,j} {\bar {b}^{ij}y_{i}^{k} y_{j}^{l}}  \right)_{(t,x)}.
\]
\noindent
Taking into account that the quadratic form  $B^{kl}$ is
nondegenerated, we obtain that $B^{11} \neq 0$ everywhere. From the
equality
\[
\dfrac{{B^{kl}}}{{B^{11}}}\left( {\tau ,y,v} \right) =
\dfrac{\sum_{i,j} {\bar {b}^{ij}y_{i}^{k} y_{j}^{l}} }
   {\sum_{i,j} {\bar {b}^{ij}y_{i}^{1}y_{j}^{1}} }\left( {t,x} \right)
\]
\noindent
we receive that this fraction is
function of $\left( {t,y} \right)$.
Thus
\[
B^{kl}(\tau ,y,v) = A(\tau ,y,v)\bar {B}^{kl}\left( {\tau ,y} \right)
\]
\noindent
for $A\left( {\tau ,y,v}
\right) = B^{11}\left( {\tau ,y,v} \right)$ and some functions $\bar
{B}^{kl}(t,y)$. Therefore, $\B  \in \Ob _{\PE_{2}}$. \qed

%\subsection*{}
\bigskip\noindent
2. $\PE_{3} =\PE_{1} \cap\PE_{2} $ is
closed in $\PE$, in $\PE_{1} $, and in
$\PE_{2} $, because $\PE_{1} $ and
$\PE_{2}$ are closed in $\PE$.  \qed

%\subsection*{}
\bigskip\noindent
3. Suppose $\A  \in\Ob _{\PE_{4}} $ and $F\colon
\A  \to \B $ is a morphism of $\PE$. From the
first equation of \eqref{eq_PE} we obtain that $B^{kl}\left( {\tau
,y,v} \right)$ is independent of $v$. Hence $B^{kl} = B^{kl}\left(
{\tau ,y} \right)$, $\PE_{4} $ is closed in $\PE$,
so it is closed in $\PE_2$ too.  \qed

%\subsection*{}
\bigskip\noindent
4. Since $\PE_{3} $ and $\PE_{4}$ are closed in $\PE$, we obtain that
$\PE_{5} =\PE_{3} \cap\PE_{4} $ is closed in $\PE$, $\PE_{3} $ and
$\PE_{4} $.
\qed
%\end{proof}

\subsection*{Proof of Theorem \ref{thm_TPE}}

%\subsection*{}
\noindent
1. By definition, $\TPE $ is wide in $\PE$.

\noindent Suppose $F\colon \A  \to \B $ is a morphism
in $\PE$. By Theorem \ref{thm_mor_PE}, the function $\tau
(t)$ is nondegenerated, so we could consider the inverse function
$t\left( {\tau}  \right)$. The map $\left( {\tau ,y,v} \right) \to
\left( {t\left( {\tau}  \right),y,v} \right)$ is an isomorphism in
$\PE$. Note that superposition of $F$ with this isomorphism
is a morphism in $\TPE $. Therefore $\TPE $ is
plentiful in $\PE$.  \qed

%\subsection*{}
\bigskip\noindent
2. $\TPE _k$ is closed in $\PE$, and
$\TPE $ is wide and plentiful in $\PE$.
Thus $\TPE _{k} =\PE_{k} \cap\TPE $ is closed in 
$\TPE $ and also it is wide and plentiful in $\PE_k$. 
\qed
%\end{proof}

\subsection*{Proof of Theorem \ref{thm_QPE}}

%\subsection*{}

Using the system \eqref{eq_PE}, we obtain that the map 
$(t,x,u) \to \left( {t,y,\varphi u + \psi}  \right)$ is a morphism 
in $\QPE$ if and only if
\begin{equation}
\label{eq_QPE}
\left\{
\begin{aligned}
 B^{kl} & = \sum_{i,j} b^{ij}y_{i}^{k} y_{j}^{l}  \\
 B^{k} & = \sum_{i,j} b^{ij}y_{ij}^{k} 
    + 2 \sum_{i,j} b^{ij}\left( \ln \bar {\varphi} \right)_j y_i^k
   + \sum_{i} b^{i}y_{i}^{k} - y_{t}^{k}  \\
 Q\bar {\varphi} & = 
  \left( \sum_{i,j} {b^{ij}\bar {\varphi} _{ij} 
    + \sum_{i} b^{i} \bar{\varphi}_i - \bar {\varphi}_t} \right)v
   + \left( \sum_{i,j} b^{ij}\bar{\psi}_{ij} + \sum_{i} b^i\bar {\psi}_i 
   - \bar {\psi} _{t}  \right)
   + q\left( {t,x,\bar {\varphi} v + \bar {\psi} } \right)
\end{aligned}
\right.,
\end{equation}
\noindent where $\bar {\varphi} = \varphi ^{-1}$, $\bar {\psi} = -
\varphi ^{-1}\psi$, so $U = \bar {\varphi} v + \bar {\psi} $. By
definition, all subcategories of $\overline {{\QPE}}$
considered in the Theorem are full subcategories of $\overline
{{\QPE}}$.

%\subsection*{}
\bigskip\noindent
1a. If $c^{ij} = 0$ and $v$ is linear in $u$, then $C^{kl} = 0$.
By the second equation of the system \eqref{eq_PE}, it follows that
$\QPE$ is closed in $\overline {{\QPE}}$.

%\subsection*{}
\bigskip\noindent
1b. Let $F\colon \A  \to \B $, $(t,x,u)
\mapsto \left( {t,y\left( {t,x} \right),v(t,x,u)} \right)$ be a
morphism in $\TPE _1$, and $\A , \B  \in \Ob
_{\QPE}$. Using the second equation of the system
\eqref{eq_PE}, we get $\left( {\ln  U_{v}} \right)_{v} B^{kl} =
C^{kl} = 0$. It follows that $U$ is linear in $v$, $v$ is linear in
$u$, $F$ is a morphism in $\QPE$, and $\QPE$ is
full in $\TPE $.

%\subsection*{}
\bigskip\noindent
1c. Suppose $\A  \in \Ob _{\TPE _{1}}$. Fix
$u_{0} \in \Omega _{\A } $ and consider the map $F\colon
(t,x,u) \mapsto \left( {t,x,v(t,x,u)} \right)$,
\[
v(t,x,u) = \int\limits_{u_{0}} ^{u} {\exp\left(
{\int\limits_{u_{0}} ^{\xi}  {\lambda \left( {t,x,\varsigma}
\right)d\varsigma} }  \right)d\xi} .
\]
\noindent
$F$ define the isomorphism in $\TPE _1$ from
$\A $ to $\B $ with
\[
C^{ij} = \left( {\ln U_{v}}
\right)_{v} b^{ij} + U_{v} \lambda b^{ij} = v_{u}^{ - 1} \left(
{\lambda - \left( {\ln v_{u}} \right)_{u}} \right) = 0.
\]
\noindent
Therefore
every object of $\TPE _{1} $ is isomorphic in
$\TPE _{1} $ to some object of $\QPE$, and
$\QPE$ is full in $\TPE _{1}$.  \qed

%\subsection*{}
\bigskip\noindent
2. The image of a compact under a continuous map is compact. The
surjectivity of the map completes the proof.  \qed

%\subsection*{}
\bigskip\noindent
3. $\TPE _{3} $ is closed in $\mathcal{�PE}_{1} $,
$\QPE$ is fully dense in $\mathcal{�PE}_{1} $.  \qed

%\subsection*{}
\bigskip\noindent
4. $\TPE _{5} $ is closed in $\TPE _{3} $,
and $\QPE'$ is fully dense in $\TPE _{3} $.
Equality $\QPE'_{1} =\TPE _{5} \cap\QPE'$ concludes the proof. 

%\subsection*{}
\bigskip\noindent
5. Let $\A  \in \Ob _{\QPE''}$, and $F\colon \A  \to \B $ be 
a morphism in $\QPE'$. From the first equation of the system 
\eqref{eq_QPE} we obtain
\begin{equation}
\label{eq_aA}
a(t,x,u) = A(t,y,v)\bar {a}(t,x),
\end{equation}
\noindent
where 
$\bar {a}(t,x) = {{B^{11}\left( {t,y\left( {t,x}
\right)} \right)} \left/ 
{{\left( {\sum_{i,j} b^{ij}(t,x)y_i^1 y_j^1 (t,x)} \right)}} \right.}$.

From the second equation of \eqref{eq_QPE} we obtain
\begin{equation}
\label{eq_omega_k}
B^{k}(t,y,v) = A(t,y,v)\omega ^{k}\left( {t,x}
\right) + \mu ^{k}(t,x),
\end{equation}
\noindent
where
\[
\omega ^{k}(t,x) = 
\bar {a}\left( {\sum_{i,j} \bar {b}^{ij} y_{ij}^{k} 
  + 2 \sum_{i,j} \bar {b}^{ij} \left( {\ln \bar {\varphi}}  \right)_{j} y_{i}^{k} 
  + \sum_{i} \bar {b}^{i}y_{i}^{k}}  \right),
\quad
\mu ^{k}(t,x) = \sum_{i} \xi ^{i}y_{i}^{k} - y_{t}^{k}.
\]
\noindent
Further we need the following statement.

\begin{figure}[tbh]
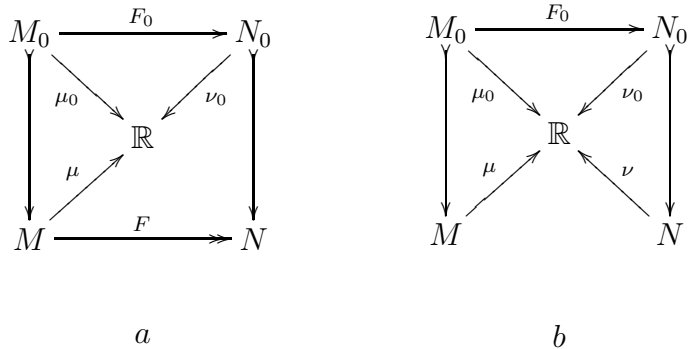

\begin{center}
$\diagram
M_0 \ddto |<<\tip \drto_{\mu_0} \rrto^{F_0} & &
N_0 \ddto |<<\tip \dlto^{\nu_0} \\
& \R & \\
M \urto^{\mu} \rrto^{F}|>>\tip  & &
N \\
& a &
\enddiagram$
\quad \quad \quad \quad
$\diagram
M_0 \ddto |<<\tip \drto_{\mu_0} \rrto^{F_0} & &
N_0 \ddto |<<\tip \dlto^{\nu_0} \\
& \R & \\
M \urto^{\mu} & &
N \ulto_{\nu} \\
& b &
\enddiagram$
\end{center}
\caption{The extension of a function}
\label{FigExt}
\end{figure}

\begin{lem}[about the extension of a function]\label{lem_lift_fun}

Suppose $M$, $N$ are $C^r$-manifolds, $1 \leq r \leq \infty$,
$F\colon M \to N$ is a surjective $C^r$-submersion, $\mu \colon M
\to \R $ is a $C^s$-function, $0 \leq s \leq r$ (if $s=0$ then $\mu$
is continuous). Take
\[
N_0 = \left\{ n \in N\colon  \;\left. {\mu} \right|_{F^{-1}(n)} =
 \const \right\},
\]
\[
M_0 = F^{-1} \left( {N_0} \right) =
\left\{ {m \in M\colon  \;\forall m' \in M \;
\left[ {F\left( {m'} \right) = F\left( m \right)} \right] \Rightarrow
\left[ {\mu \left( {m'} \right) = \mu \left( m \right)} \right]} 
\right\},
\]
\noindent
$F_0 = \left. {F} \right|_{M_0}$,
$\mu_0 = \left. {\mu} \right|_{M_0}$,
and define the function $\nu_0 \colon  N_0 \to \R$ by the formula
$\nu_0 F_0 = \mu_0$ (see Fig.~\ref{FigExt}(a)).
Then $\nu_0$ can be extended from $N_0$ to the entire manifold $N$
so that the extended function $\nu \colon  N \to \R$
will have class $C^s$ of smoothness
(see Fig.~\ref{FigExt}(b); 
both diagrams Fig.~\ref{FigExt}(a, b) are commutative).
\end{lem}

\subsection*{Proof of Lemma \ref{lem_lift_fun}}

Take an open covering $\left\{ {V_{i} \colon i \in I} \right\}$
of $N$ such that for every $V_{i}$ there exist a $C^r$-smooth section 
$p_{i} \colon V_{i} \to M$ over $V_{i}$, $F \circ p_{i} = \left.
{\id} \right|_{V_{i}}$ (such a covering exists, because $F$ is
submersive and surjective). There exist a $C^r$-partition of unity
$\left\{ {\lambda _{i}} \right\}$ that is subordinated to $\left\{
{V_{i} } \right\}$ \cite{Hirsch}. Define the functions
\[
\nu _{i} \left( {n} \right) = \left\{
{{\begin{array}{*{20}l}
 {\lambda _{i} \left( {n} \right)\mu \left( {p_{i} \left( {n} \right)}
\right),} & {n \in V_{i}}  \\
 {0,} & {n \notin V_{i}}  \\
\end{array}} } \right..
\]
\noindent
Then $\nu \left( {n} \right) = \sum\limits_{i \in I} {\nu _{i}
\left( {n} \right)} $ is the desired function.  \qed

%\end{proof}

\subsection*{Proof of Theorem \ref{thm_QPE} (continuation)}

%\subsection*{}
Fix $k$. In the notations and assumption of Lemma
\ref{lem_lift_fun}, let us replace the map $F$ by $(t,x) \mapsto
\left( {t,y(t,x)} \right)$  and the continuous function $\mu$ by
$\mu ^{k}(t,x)$. We obtain that there exists a continuous function
$\nu ^{k}\left( {t,y} \right)$ such that for each $\left( t_0,y_0
\right)$ if $\mu ^{k}(t,x)$ is constant on the pre-image of $\left(
t_0,y_0 \right)$ with respect to the map $(t,x) \mapsto \left(
{t,y(t,x)} \right)$ then $\nu ^{k}\left( {t_0,y_0} \right)$
coincides with this constant. Denote
\begin{equation}
\label{eq_Bk}
\bar {B}^{k}(t,y,v) = {{\left( {B^{k}(t,y,v) -
\nu ^{k}(t,y)} \right)} \left/ {{A(t,y,v)}}
\right.}.
\end{equation}

Let us consider the following two cases for every point $\left(
{t_{0} ,y_{0}} \right) $.

%\subsection*{}
\bigskip\noindent
\underline{Case 1}: $A\left( {t_{0} ,y_{0} ,v} \right)$ is
independent of $v$. Using \eqref{eq_omega_k}, we obtain that
$B^{k}\left( {t_{0} ,y_{0} ,v} \right)$ is independent of $v$; and
using \eqref{eq_Bk} that $\bar {B}^{k}$ is independent of $v$.

%\subsection*{}
\bigskip\noindent
\underline{Case 2}: For given $\left( {t_{0} ,y_{0}}  \right)$ the
set $\left \{ A\left( {t_{0} ,y_{0} ,v} \right)\colon  v \in \Omega
\right\}$ contains more then one element. Using \eqref{eq_omega_k},
we obtain that the restriction of $\mu ^{k}\left( {t_{0} ,x}
\right)$ to the pre-image of the point $\left( {t_{0} ,y_{0}}
\right)$ is constant. Then $\mu ^{k}\left( {t_{0} ,x} \right) = \nu
^{k}\left( {t_{0} ,y_{0}}  \right)$ on this pre-image, and $\bar
{B}^{k} = \omega ^{k}(t,x)$ is independent of $v$ in this case too.

Therefore, $B^{k}(t,y,v) = A(t,y,v)\bar {B}^{k}(t,y) + \nu^{k}(t,y)$. 
So, the equation $\B $ is of the form
\[
v_{t} = A(t,y,v)\left( {\sum_{k,l}\bar {B}^{kl}(t,y)v_{kl}
 + \sum_{k}\bar {B}^{k}(t,y)v_{k}}  \right) 
 + \sum_{k}\nu^{k}(t,y)v_{k} + Q(t,y,v),
\]
\noindent and $\B $ is the object of $\QPE''$.

For $F$ to be a morphism in $\QPE''$, it is necessary and
sufficient to have
\begin{equation}
\label{eq_QPE''}
\left\{
\begin{aligned}
& a(t,x,u) = A(t,y,v)\bar {a}(t,x)
\\
& \bar {B}^{kl} (t,y)
  =  \bar {a} \sum_{i,j}\bar {b}^{ij}y_{i}^{k} y_{j}^{l} (t,x)
\\
&
y_t^k + \Xi^k - \sum_{i} \xi^i y_i^k 
   = a(t,x,u) \left(
   \sum_{i,j} \bar{b}^{ij}y_{ij}^{k}
   + 2\sum_{i,j} \bar{b}^{ij}\left( {\ln \bar{\varphi}} \right)_j y_i^k
   + \sum_{i} \bar {b}^{i}y_{i}^{k} -B^{k}/ \bar {a} \right)
%   (t,x)
\\
& Q\bar {\varphi} = 
  \left( {\sum_{i,j} a\bar {b}^{ij}\bar {\varphi} _{ij}
    + \sum_{i} \left( a\bar {b}^{i}+\xi_i \right) \bar{\varphi}_{i}
   - \bar {\varphi} _{t}} \right)v +
 \\
 &
 \quad \quad \quad
 + \left( {\sum_{i,j} a\bar{b}^{ij}\bar {\psi }_{ij}
  + \sum_{i} \left( a\bar {b}^{i}+\xi_i \right) \bar {\psi} _{i} 
  - \bar {\psi} _{t}}  \right) +
  q\left( {t,x,\bar {\varphi} v + \bar {\psi} } \right)
\end{aligned}
\right.
\end{equation}

%\subsection*{}
\bigskip\noindent
6. $\QPE''_1$ is closed in $\QPE''$ and in
$\QPE'_1$, because $\QPE''$ and
$\QPE'_{1} $ are closed in $\QPE'$.
$\QPE''_1$ is closed in $\QPE''_0$, because
$\QPE''_0$ is the subcategory of $\QPE''$.  \qed

%\subsection*{}
\bigskip\noindent
7. Suppose $\A  \in \Ob _{\QPE''_{1q}}$, $F\colon
\A  \to \B $ is a morphism in $\QPE''_{1} $.
From the third equation of \eqref{eq_QPE} we get
\begin{gather}
Q(t,y,v)  =
 \notag
\\ 
= \left( 
    \sum_{i,j} b^{ij}\bar{\varphi}_{ij}
     + \sum_{i} b^{i}\bar{\varphi}_{i} 
     + q_{1} (t,x) - \bar{\varphi}_{t} 
     \right) \bar{\varphi}^{-1} v
  + \left( 
    \sum_{i,j} b^{ij}\bar{\psi}_{ij}
    + \sum_{i} b^{i}\bar{\psi}_{i} + q_{0} (t,x)
  - \bar {\psi}_{t} \right) \bar{\varphi}^{-1} =
 \notag
\\ %\qquad
  =  Q_{1}(t,x)v + Q_{0}(t,x), \notag
\end{gather}
\noindent so $Q_1$, $Q_0$ are functions of $(t, y)$, and $\B
\in\Ob _{\QPE''_{1q}}$. Thus $\QPE''_{1q} $ is
closed in $\QPE''_{1} $.  \qed

%\subsection*{}
\bigskip\noindent
8. Suppose $\A  \in \Ob _{\QPE'_{n}}$, $F\colon
\A  \to \B $ is a morphism in $\QPE'$. For
given $\left( {t_{0} ,y_{0}}  \right)$ let us fix arbitrary $x_{0}$
such that $y\left( {t_{0} ,x_{0}} \right) = y_{0} $. Using $\bar
{\varphi}  \ne 0$ and $a \in \mathcal {A}_{nc} \left( {T \times X}
\right)$, from \eqref{eq_aA} we get
\[
A\left( {t_{0} ,y_{0} ,v}
\right) = a\left( {t_0 ,x_0 ,\bar {\varphi} \left( {t_0 ,x_0}
\right)v + \bar {\psi} \left( {t_0 ,x_0}  \right) } \right) \bar
{a}\left( {t_0 ,x_0}  \right) \ne \mathrm{const}.
\]
\noindent
Finally, we obtain 
$A \in \mathcal {A}_{nc} \left( {T \times Y} \right)$, and
$\B  \in\Ob _{\QPE'_{n}}$, so $\QPE'_{n}$ is closed in $\QPE'$.  \qed

%\subsection*{}
\bigskip\noindent
9. Suppose $\A  \in \Ob _{\QPE''_{0n}}$, $\B
\in \Ob _{\QPE''_{n}}$. Substituting $\xi_i=0$ in the third
equation of \eqref{eq_QPE''}, we get
\[
y_{t}^{k} + \Xi ^{k} (t,y) =
 a(t,x,u) \left(
   \sum_{i,j} \bar {b}^{ij}y_{ij}^{k}
   + 2 \sum_{i,j} \bar {b}^{ij}\left( {\ln  \bar {\varphi}} \right)_{j} y_{i}^{k}
   + \sum_{i} \bar {b}^{i}y_{i}^{k} -B^{k}/ \bar {a} \right) (t,x).
\]
\noindent Since $a \in \mathcal {A}_{nc} (T \times X)$, and left
hand side is independent of $u$, it follows that both sides of this
equality vanishes, and
\begin{equation}
\label{eq_y_t}
y_{t}^{k} = - \Xi ^{k} (t,y)
\end{equation}
The function $y(t,x)$ satisfies the ordinary differential equation
\eqref{eq_y_t} with smooth right hand side, so for any $t,t'$ an
equality $y(t,x_1)= y(t,x_2)$ implies that $y(t',x_1)= y(t',x_2)$.
Let 1-parameter transformation group $g_s\colon  T \times Y \to T \times
Y$ be given by $(t,y(t,x)) \mapsto (t+s,y(t+s,x))$. This group is
correctly defined when $T=\R$; otherwise transformations $g_s$ are
partially defined, nevertheless reasoning below remains correct
after small refinement.

For every $s$ the composition $g_s g_{-s}$ is identity, so $g_s$ is
bijective. $\{g_s\}$ is a flow map of the smooth vector field
$\partial _t - \sum_{k} \Xi^{k}(t,y) \partial _{y^k}$, 
so transformations $\{g_s\}$ are smooth by both $t$ and $y$.

Define the map $z(t,y)$ by the equality $g_{-t}(t,y)=(0,z(t,y))$.
Then the map $G\colon  T \times Y \to T \times Y$, $(t,y) \mapsto
(t,z(t,y))$ is the isomorphism in $\QPE''$ such that for
every $x,t$ $z(t,y(t,x))=z(0,y(0,x))$. Therefore $G \circ F \in \Hom
_{\QPE''_{0}}$.  \qed

%\subsection*{}
\bigskip\noindent
10. Let $\A$ be an object of $\QPE_c''$. 
Consider the solution $y\colon T \times X \to X$ of the linear PDE 
${{\partial y^{k}}\left/ {{\partial t}} \right.} = 
 \sum_{i} \xi ^{i}(t,x){{\partial y^{k}} \left/ {{\partial x^{i}}} \right.}$ 
(the solution exists in light of the compactness of $X$). 
The isomorphism $(t,x,u) \mapsto \left( {t,y(t,x),u} \right)$ maps $\A$ 
to some object of $\QPE''_{0} $. 
Thus $\QPE''_{0c} $ is closed in $\QPE''_c$. 
\qed
%\end{proof}

\subsection*{Proof of Theorem \ref{thm_QPE_a}}

%\subsection*{}

If $a \neq \mathrm{const}$, then $\QPE''_{0a} (a)$ is fully
plentiful in $\QPE''_{a} (a)$ thanks to the part 9 of
Theorem \ref{thm_QPE}.

If $a = \mathrm{const}$, then $\QPE''_{a}(a)$ coincides with
$\QPE''_{1} $,
which is closed in $\QPE''$ by the Theorem~\ref{thm_QPE}.
So $\QPE''_{a} (a)$ is fully plentiful
in $\QPE''$.

Suppose now that $a \ne \mathrm{const}$, $\A  \in \Ob
_{\QPE''_{a} (a)}$, and $F\colon \A  \to \B $
is a morphism in $\QPE''$. Consider the equation
\eqref{eq_aA} as functional one:
\begin{equation}
\label{eq_a_func}
a\left( {\bar {\varphi} (t,x)v + \bar {\psi} \left( {t,x}
\right)} \right) = A(t,y,v) \bar {a} (t,x).
\end{equation}
Let us consider the following three cases.

%\subsection*{}
\bigskip\noindent
\underline{Case 1.} $a(u) = He^{\lambda u}$, $\lambda , H =
\mathrm{const}$, $\lambda \ne 0$. Substituting the formula for $a$
to \eqref{eq_a_func}, we get $\lambda \bar {\varphi} (t,x)v - \ln
A(t,y,v) = \left( {\ln \bar {a} - \lambda \bar {\psi}  - \ln H}
\right)$. The right hand side of the equality is a function of
$(t,x)$, so $\bar {\varphi}  = \bar {\varphi} (t,y)$, and the
isomorphism $(t,y,v) \mapsto \left( {t,y,\bar {\varphi} (t,y)v}
\right)$ maps $\B $ to some object of $\QPE''_{a}
(a)$.

%\subsection*{}
\bigskip\noindent
\underline{Case 2.} $a(u) = H\left( {u - u_{0}}  \right)^{\lambda}
$, $\lambda ,H,u_{0} = \mathrm{const}$, $\lambda \ne 0$.
Substituting the formula for $a$ to \eqref{eq_a_func}, we get
\[
\left( {v + \bar {\varphi}^{-1}(t,x)
\left( {\bar {\psi} (t,x) - u_{0}}  \right)} \right)^{\lambda}
 = A(t,y,v)H^{ - 1}\bar {\varphi} ^{ - \lambda} \bar {a}(t,x).
\]
\noindent Thus $\bar{\varphi}^{-1}\left({\bar{\psi}-u_0} \right) =
q(t,y)$ for some function $q$, and the object $\B $ maps by the
isomorphism $\left( {t,y,v} \right) \mapsto \left( {t,y,v + q(t,y) +
u_{0}} \right)$ to some object of $\QPE''_{a} (a)$.

%\subsection*{}
\bigskip\noindent
\underline{Case 3.} Suppose now that $a(u)$ is neither $He^{\lambda
u}$ nor $H\left( {u - u_{0}}  \right)^{\lambda} $. Denote $\bar {x}
= (t,x)$, $\bar {y} = (t,y)$, $\alpha = \ln a$. Fix a point $\bar
{y}_{0} $ and take $Z = \left\{ {\bar {x}\colon \bar {y}\left( {\bar {x}}
\right) = \bar {y}_{0}} \right\} \subset T \times X$. Using
\eqref{eq_a_func}, we obtain that $\forall \bar {x}_{0} ,\bar
{x}_{1} \in Z\quad \alpha \left( {\bar {\varphi} _{1} z + \bar
{\psi} _{1}} \right) - \alpha \left( {\bar {\varphi} _{0} z + \bar
{\psi} _{0}} \right)$ is independent of $v$, where $\bar {\varphi}
_{i} = \bar {\varphi} \left( {\bar {x}_{i}} \right)$, $\bar {\psi}
_{i} = \bar {\psi} \left( {\bar {x}_{i}} \right)$). Consider
additive subgroup $G = G\left( {\bar {y}_{0}} \right)$ of $\R$
generated by the set $\left\{ {\ln \bar {\varphi }\left( {\bar {x}}
\right) - \ln \bar {\varphi} \left( {\bar {x}_{0}} \right)\colon \bar {x}
\in Z} \right\}$.

Consider the following two subcases.

%\subsection*{}
\bigskip\noindent
\underline{Case 3.1.} $G \ne \left\{ {0} \right\}$.

\noindent
Put $\hat H_{1} = \ln\bar{\varphi}_1 - \ln\bar{\varphi}_0 \in G - \{0\}$, 
$u_{0} = {{\left( {\bar {\psi} _{0} - \bar {\psi} _{1}} \right)} \left/ 
{{\left( {\bar {\varphi} _{1} - \bar {\varphi} _{0}} \right)}} \right.}$. 
Substituting $v = {{\left({w + u_0 - \bar {\psi}_0} \right)} \left/ 
{{\bar {\varphi} _{0}} } \right.}$, for any $w$ we have 
$\alpha \left( {e^{\hat H_1} w + u_0}\right) - 
\alpha \left( {w+u_0}\right) = c = \const$. 
Consider the function 
$\beta \left( {x} \right) = \alpha \left( {e^{x} + u_{0}}\right)$. 
Using $\beta \left( {x + \hat H_{1} } \right) = \beta (x)+ c$, 
we obtain that the function $\beta \left( {x} \right) - \lambda x$ is  
$\hat H_{1}$-periodic, where $\lambda = {c}/{\hat H_{1}}$. Then
\[
a(u) = \left( {u - u_{0}}\right)^{\lambda}
H\left( {\ln \left( {u - u_{0}}  \right)}\right),
\]
\noindent
where $H$ is $\hat H_{1}$-periodic, $H \ne \mathrm{const}$,
because case ``$H = \mathrm{const}$'' was already considered. 
Let $\hat H > 0$ be the smallest positive period of $H$.
For all $\bar {x} \in Z$ the number $\ln \bar {\varphi} \left( {\bar
{x}} \right) - \ln \bar {\varphi} _{0} $ is a multiple of $\hat H$,
so for any $\bar {y}_{0}$ we have
$\bar {\varphi} \left( {\bar {x}} \right)
\in \left\{ {\bar {\varphi} _{0} \exp \left({k \hat H}\right)\colon
k \in \mathbb{Z}}\right\}$.
Note that $\hat H$ is independent of $\bar {y}_{0}$,
because $a(u)$ is independent of $\bar {y}_{0} $.

%\subsection*{}
\bigskip\noindent
\underline{Case 3.2.} $G = \left\{ {0} \right\}$, that is $\left.
{\bar {\varphi} } \right|_{Z} \equiv \bar {\varphi} _{0} =
\mathrm{const}$. Now we have the following two subsubcases:

%\subsection*{}
\bigskip\noindent
\underline{Case 3.2.a.} $\left. {\bar {\psi} } \right|_{Z} \ne
\mathrm{const}$, that is $\exists \bar {x}_{0} ,\bar {x}_{1} \in Z\colon
\bar {\psi} \left( {\bar {x}_{1}} \right) - \bar{\psi} \left({\bar
{x}_{0}} \right) = \hat H_{1} \ne 0$. Then $\alpha \left( {u + \hat
H_{1}}  \right) - \alpha (u) = \mathrm{const}$. By the same token as
in case 3.1 we get $a(u) = H(u)e^{\lambda u}$, where $\lambda =
\mathrm{const}$, $H$ is a periodic function with the smallest period
$\hat H > 0$. Note that such representation of $a(u)$ is unique.
Substituting this to \eqref{eq_a_func}, we obtain that $\forall \bar
{y}\;\forall \bar {x}_{0} ,\bar {x}_{1} \in Z_{\bar {y}}$ the number
$\bar {\psi} \left( {\bar {x}_{1}}  \right) - \bar {\psi} \left(
{\bar {x}_{0}}  \right)$ is a multiple of $\hat H$.

%\subsection*{}
\bigskip\noindent
\underline{Case 3.2.b.} $\left. {\bar {\psi} } \right|_{Z} =
\mathrm{const}$ for given $\bar {y}_{0}$. We already considered thes
cases $a(u) = H(u)e^{\lambda u}$ and $a(u) = \left( {u - u_{0}}
\right)^{\lambda} H\left( {\ln \left({u-u_{0}} \right)} \right)$, so
we can assume without loss of generality that $a$ is not of this
form. Then at every $\bar {y}_{0}$ we have $\left. {\bar {\psi} }
\right|_{Z} = \mathrm{const}$, $\bar {\varphi}  = \bar {\varphi}
\left( {\bar {y}} \right)$, and $\bar {\psi}  = \bar {\psi} \left(
{\bar {y}} \right)$. Thus the isomorphism $(t,y,v) \to \left(
{t,y,\bar {\varphi} (t,y)v + \bar {\psi} \left( {t,y} \right)}
\right)$ maps $\B $ to some object of $\QPE''_{a} \left( {a}
\right)$.

The proof of the full density of $\QPE''_{0ca} (a)$ in
$\QPE''_{ca} \left( {a} \right)$ is similar to the proof of
part 10 in Theorem \ref{thm_QPE}.  \qed

%\end{proof}

\subsection*{Proof of Theorem \ref{thm_SQPE}}

%\subsection*{}
\bigskip\noindent
1. $\QPE''$ is closed in $\overline {{\QPE}} $, and
$\overline {{\SQPE}} $ is the subcategory of $\overline
{{\QPE}} $. Therefore $\SQPE$ is closed in
$\overline {{\SQPE}} $.  \qed

%\subsection*{}
\bigskip\noindent
2. $\SQPE_{n} $ is closed in $\overline {{\SQPE}} $
for the same reason as in Part 1 of this Theorem. This implies that
$\SQPE_{n} $ is closed in $\SQPE$.

Suppose $\A$ is an object of $\SQPE_{0}$, 
$F\colon \A  \to \B $ is a morphism in $\SQPE$. Then
$B^k(t,y,v) = A(t,y,v)\omega ^{k}(t,x)$, where $\omega ^{k}$ is
defined as in \eqref{eq_omega_k}. Hence $\omega ^{k}$ is a function
of $(t,y)$, and $B$ is a object of $\SQPE_0$.

Suppose $\A$ is an object of $\SQPE_{b}$, $F\colon \A  \to \B $  
is a morphism in $\SQPE$. From the first equation of \eqref{eq_QPE} 
we obtain that
\[\dfrac{{\bar
{B}^{kl}}}{{\bar {B}^{11}}}(t,y) = 
  \dfrac{{\sum_{i,j} \bar {b}^{ij}y_{i}^{k} y_{j}^{l} }}
    {{\sum_{i,j} \bar {b}^{ij}y_{i}^{1} y_{j}^{1}} }(x).
\]
\noindent
The right hand side is independent of $t$, so it is a function of $y$;
denote this function by ${\bar {B}}'^{kl}(y)$. Then $A\bar {B}^{kl} =
{A}'(t,y,v){\bar {B}}'^{kl}(y)$, where ${A}' = AB^{11}$. 
It follows that $\B$ is an object of $\SQPE_{b}$, 
and $\SQPE_{b} $ is closed in $\SQPE$.  \qed

%\subsection*{}
\bigskip\noindent
3. Let us recall that $\SQPE_{0n}$ is closed in
$\QPE''_{0n} $. So it suffices to prove that any morphism
in $\QPE''_{0n} $ is also a morphism in
$\SQPE_{0n}$. Suppose $F\colon \A  \to \B $
is a morphism in $\QPE''_{0n} $. Then $y_{t}^{k} \left(
{t,x} \right) = A(t,y,v)\omega ^{k}(t,x)$, where
\[
\omega ^{k} = -
\bar {B}^{k} + \bar{a} \left( {
  \sum_{i,j} \bar {b}^{ij}y_{ij}^{k} 
  + 2 \sum_{i,j} \bar{b}^{ij}\left( {\ln\bar{\varphi}} \right)_{j} y_{i}^{k} 
  + \sum_{i,j} \bar{b}^{i}y_{i}^{k}} \right).
\]
\noindent
Since the left hand side of this
equality is independent of $v$ and $A \in \mathcal {A}_{nc} (Y)$,
we conclude that $\omega ^{k}=0$. Thus $F$ is a morphism in
$\SQPE_{0n} $. Finally, $\SQPE_{0n} =
\QPE''_{0n} $, is closed in $\QPE''_{0} $ and is
fully dense in $\QPE''_{n} $.  \qed

%\subsection*{}
\bigskip\noindent
4. $\QPE''_{1} $ is closed in $\overline {{\QPE}} $, so
$\SQPE_{1} $ is closed in $\overline {{\SQPE}} $ and, consequently,
is closed in $\SQPE_{0} $.  \qed

%\subsection*{}
\bigskip\noindent
5. The proof is similar to the proof of part 1 of Theorem \ref{thm_QPE_a}.  
\qed

%\end{proof}

\subsection*{Proof of Theorem \ref{thm_AQPE}}

%\subsection*{}

From \eqref{eq_QPE}-\eqref{eq_aA} and the fact that
$\SQPE_{b} $ is closed in $\overline {{\SQPE}} $ it
follows that the map $\left( {t,x,u} \right) \mapsto \left(
{t,y,\varphi u + \psi}  \right)$ is a morphism in $\SQPE$
with the source from $\AQPE$ if and only if the following
conditions are satisfied:
\begin{equation}
\label{eq_AQPE}
\left\{
\begin{aligned}
 & A(t,y,v) & = & a(x,u)\bar{a}(t,x)
\\
 & \bar {B}^{kl}(y) & = & \bar {a}(t,x)\nabla y^{k}\nabla y^{l}
\\
 & B^{k}(t,y,v) & = & A(t,y,v)\bar {B}^{k} (t,y) + C^{k}(t,y) =
\\
& &  = & a(x,u)\left( {\Delta y^k+\left({\eta
  + 2\nabla \left({\ln\bar{\varphi}}\right)}\right)\nabla y^{k}} \right)
 + \xi \nabla y^{k}
\\
 & Q\bar \varphi & = & \left( {a
 \left( {\Delta \bar \varphi + \eta \nabla \bar \varphi } \right)
 + \xi \nabla \bar \varphi - \bar \varphi _t } \right)v +
\\
 & & & + \left( {a\left( {\Delta \bar \psi + \eta \nabla \bar \psi } \right)
 + \xi \nabla \bar \psi - \bar \psi _t } \right)
 + q\left( {t,x,\bar \varphi v + \bar \psi } \right)
\end{aligned}
\right.
\end{equation}

%\subsection*{}
\bigskip\noindent
1. Suppose $F\colon \A  \to \B $ is a morphism in $\overline {\AQPE}$, 
$\A$ is an object of $\AQPE$. From the second equation of the system
\eqref{eq_AQPE} it follows that $\bar {a} = \bar {a}(x)$. Using the
first equation of \eqref{eq_AQPE} and taking into account the
independence of $\bar {\varphi}$, $\bar {\psi} $ on $t$, we get the
independence of $A = A\left( {y,v} \right)$ on $t$. It follows from
the third equation of \eqref{eq_AQPE} that $B^{k}$ is independent of
$t$, and $B^{k}(y,v) = A(y,v)\bar {B}^{k}(t,y) + C^{k}\left( {t,y}
\right)$. From this formula, by the same token as in proof of part 4
of Theorem \ref{thm_QPE} we obtain existing of functions ${\rm
H}^{k}(y)$, $\Xi ^{k}(y)$ such that $B^{k} = A(y,v){\rm H}^{k}(y) +
\Xi ^{k}(y)$. Substituting $u = \bar {\varphi} (x)v + \bar {\psi}
(x)$ in the last equation of \eqref{eq_AQPE}, we obtain that $Q$ is
independent of $t$. This implies that target $\B$ of the morphism $F$ is
of the form
\[
v_{t} = A(y,v)\left( { \sum_{k,l} \bar {B}^{kl}(y)v_{kl} 
  + \sum_{k} {\rm H}^{k}(y)v_{k}} \right) 
  + \sum_{k} \Xi^{k}(y)v_{k} + Q(y,v).
\]
\noindent
We received this form of $\B$ only locally. Neverytheless we can lead 
it to the equation of the same form but with globally defined function
$A(y,v)$, for example by the way described in Remark \ref{rem:a}. 
Then quadratic form $\bar {B}^{kl}$ is defined on the whole manifold $Y$,
so we can equip $Y$ with Riemann metric $\bar {B}^{kl}$ and finally get 
$\B  \in \Ob _{\AQPE}$.  \qed

%\subsection*{}
\bigskip\noindent
2. $\AQPE_{n} = \AQPE \cap \SQPE_{n} $ is
closed in $\AQPE$, because $\SQPE_{n} $ is closed in
$\SQPE$.

Let $F\colon \A  \to \B $ be a morphism in
$\SQPE_{bn}$, and both source and target of $F$ are objects
of $\AQPE_{n} $. Then $\bar {a}$ is independent of $t$, and
\begin{equation}
\label{eq_AQPE_func}
a\left( {x,\bar {\varphi} (t,x)v + \bar {\psi} \left( {t,x} \right)}
 \right) =
A\left( {y(x),v} \right) \bar {a}(x).
\end{equation}
Let $x = x_{0}$. Suppose that the set $\left\{ {\left( {\bar
{\varphi }\left( {t,x_{0}}  \right), \bar {\psi} \left( {t,x_{0}}
\right)} \right)} \right\}$ have more than one element, and consider
the intervals
\[
I\left( {v} \right) = \left\{ {\left( {\bar {\varphi}
\left( {t,x_{0}} \right)v + \bar {\psi} \left( {t,x_{0}}  \right)}
\right) \colon t \in T_{\A }}  \right\} \subseteq \R.
\]
\noindent
Then
$a\left( {x_{0} ,u} \right)$ is constant on any interval  $u \in
I\left( {v} \right)$, because the right hand side of
\eqref{eq_AQPE_func} is independent of $t$. Note that $I\left( {v}
\right)$ is a continuous function of $v$ in the Hausdorff metric,
and $\forall t\;\bar {\varphi} \left( {t,x_{0}}  \right) \ne 0$. If
at any $v$ the interval $I\left( {v} \right)$ does not collapses
into a point, then $a\left( {x_{0} ,u} \right)$ is constant on
$\bigcup I\left( {v} \right)$. But this contradicts to the condition 
$a \in \mathcal {A}_{nc} \left( {X} \right)$. 
Therefore $I\left( {v_{0}} \right)$ degenerate into the point at some 
$v_{0}$: $\bar {\varphi} \left( {t,x_{0}}  \right)v_{0} + \bar {\psi} 
\left( {t,x_{0}} \right) \equiv u_{0} $, and 
$\bar {\varphi} v + \bar {\psi}  = \bar {\varphi} \left( {t,x_{0}} 
\right)\left( {v - v_{0}} \right) + u_{0}$. 
By the assumption, ${\rm card} \left\{ {\left( {\bar {\varphi}
\left( {t,x_{0}}  \right),\bar {\psi} \left( {t,x_{0}} \right)}
\right)} \right\} >1$, so the set $\left\{ {\bar {\varphi} \left(
{t,x_{0}} \right)} \right\}$ is nondegenerated interval. Therefore,
$a\left( {x_{0} ,u} \right)$ is constant on the sets $\left\{ {u <
u_{0}} \right\}$ and $\left\{ {u > u_{0}} \right\}$. But this
contradicts to the condition $a \in \mathcal {A}_{nc} (X)$ and
continuity of $a$. This contradiction shows that for each $x_{0} $
the functions $\bar {\varphi }$, $\bar {\psi} $ are independent of
$t$. Consequently $F$ is a morphism in $\AQPE$, and
$\AQPE_{n} $ is the full subcategory of
$\SQPE_{bn} $.  \qed

%\subsection*{}
\bigskip\noindent
3. $\AQPE_{0} $ and $\AQPE_{1} $ are closed in
$\AQPE$, because $\SQPE_{0} $ and
$\SQPE_{1} $ are closed in $\SQPE$.  \qed

%\subsection*{}
\bigskip\noindent
4. If $a \notin \mathcal{A}_{\exp} \cup \mathcal{A}_{\mathrm{deg}}$,
then $\AQPE_{a} (a)$ is plentiful in $\AQPE$ by
the same arguments as used in the proof of part 1 of Theorem
\ref{thm_QPE_a}, after replacement of $\bar {x}$, $\bar {y}$ to
$x,y$ respectively.  \qed

%\subsection*{}
\bigskip\noindent
5. Let $F\colon \A  \to \B $ be a morphism in $\SQPE_{na} (a)$, 
$\A$ be an object of $\AQPE_{na}(a)$. Then
\[
a\left( {\bar {\varphi} (t,x)v + \bar {\psi} (t,x)}
\right) = A(v) \bar {a}(x).
\]
\noindent As we proved in part 2, $\bar {\varphi},\bar {\psi} $ are
independent of $t$, $F$ is a morphism in $\AQPE$, and $\B  \in \Ob
_{\AQPE} \cap \Ob _{\SQPE_{na} (a)} = \Ob _{\AQPE_{na} (a)}$. Since
$\AQPE_{na} (a)$ is full in $\AQPE_{n} $, we see that $F$ is a
morphism in $\AQPE_{na} (a)$. \qed

%\end{proof}

\subsection*{Proof of Theorem \ref{thm_EPE}}

%\subsection*{}
\bigskip\noindent
1. $\EPE$ is closed in $\overline {\EPE}$, because
$\AQPE$ is closed in $\overline {\AQPE}$.  \qed

%\subsection*{}
\bigskip\noindent
2. $\EPE_{n} $, $\EPE_{0} $, $\EPE_{1}$ are closed in $\EPE$, because 
$\AQPE_{n} $, $\AQPE_{0} $, $\AQPE_{1} $ are closed in $\AQPE$.  \qed

Suppose $F\colon \A  \to \B $ is a morphism in
$\EPE$, and $\A  \in \Ob _{\EPE_{a} (a)}$.
Then the first equation of \eqref{eq_QPE} is of the form $A(y,u)
\bar {B}^{kl}(y) = a(u) \nabla y^{k} \nabla y^{l} $. Hence $\nabla
y^{k} \nabla y^{l} = g^{kl}(y)$ for some functions $g^{kl}$. For
$\bar {B}^{kl} = g^{kl}(y)$ we get $A(y,u) = a(u)$. 
So $\A$ is the object of $\EPE_{a} (a)$, and $\EPE_{a} (a)$ is
closed in $\EPE$.  \qed

%\subsection*{}
\bigskip\noindent
3. The proof is similar to the proof of Theorem \ref{thm_QPE}.
\qed
%\end{proof}

%\begin{thebibliography}{30}


\begin{thebibliography}{99}

\bibitem [Elkin, 1999]{Elkin}
V. I. Elkin (1999), 
\textit{Reduction of nonlinear control systems. A differential 
geometric approach}, 
Mathematics and its Applications, \textbf{472}, 
Kluwer Academic Publishers, Dordrecht.

\bibitem [Hirsch, 1976]{Hirsch}
M. W. Hirsch (1967), 
\textit{Differential topology}, Graduate texts in mathematics, 
\textbf{33}, Springer-Verlag, NY.

\bibitem [Krasil'shchik, 1999]{Diffeotopy}
J. Krasil'shchik, A. Vinogradov (eds.) (1999),
%V. N. Chetverikov, A. B. Bocharov, V. N. Chetverikov, S. V. Duzhin et all, 
\textit{Symmetries and conservation laws for differential equations of 
mathematical physics}, 
Translations of Mathematical Monographs, \textbf{182}, AMS, Providence, RI.

\bibitem [Mac Lane, 1998]{Maclain}
S. Mac Lane (1998), 
\textit{Categories for the working mathematician}, Springer-Verlag, NY.

\bibitem [Olver, 1993]{Olver}
P. J. Olver (1993), 
\textit{Applications of Lie groups to differential equations}, 
Springer-Verlag, NY.

\bibitem [Prokhorova, 1998]{Prokh1998}
M. F. Prokhorova (1998), 
\textit{Modeling of solutions of heat equation and Stefan problem with 
dimension decrease}, 
Russian~Acad.~Sci.~Docl.~Math. \textbf{58} (1), 88--90.

\bibitem [Prokhorova, 2000]{Prokh2000}
M. F. Prokhorova (2000), 
\textit{Modeling of heat equation and Stefan problem}, 
Preprint VINITI N~347-V00, IMM UrBr Russian Acad. Sci., Ekaterinburg, 
1--66 (in Russian).

\bibitem [Prokhorova, 2001]{Prokh2001}
M. F. Prokhorova (2001), 
\textit{Factorization of nonlinear heat equation posed on Riemann 
manifold}, 
arXiv.org:math.AP/0108001, 1--9.

\end{thebibliography}
\end{document}